\documentclass[11pt]{article}
\usepackage{amsfonts}
\usepackage{color}
\usepackage{graphics}
\usepackage{indentfirst}
\usepackage{cite}
\usepackage{latexsym}
\usepackage{amsmath}
\usepackage{amssymb}
\usepackage[dvips]{epsfig}
\usepackage{amscd}

\usepackage{leqno} 
\newtheorem{theorem}{Theorem}[section]

\newtheorem{lemma}[theorem]{Lemma}

\newtheorem{proposition}[theorem]{Proposition}

\newcommand{\no}{\nonumber\\}
\newcommand{\n}{\varrho}
\newcommand{\nn}{\varrho_0^R}

\newcommand{\ti}{\tilde}
\newcommand{\mr}{\mathbb{R}}

\def\pf{{\it Proof.}  }
\renewcommand{\div}{ {\rm div }  }

\newcommand{\uu}{\mathbf{u}}
\newcommand{\hh}{\mathbf{H}}
\newcommand{\na}{\nabla }
\newcommand{\vp}{\varphi }
\newcommand{\xx}{\mathbf{x}}

\newcommand{\bl}{\begin{lemma}}
\newcommand{\el}{\end{lemma}}

\newcommand{\ve}{\varepsilon}
\newcommand{\la}{\label}

\newcommand{\om}{\Omega}

\newcommand{\bn}{\begin{eqnarray}}
\newcommand{\en}{\end{eqnarray}}
\newcommand{\bnn}{\begin{eqnarray*}}
\newcommand{\enn}{\end{eqnarray*}}

\newcommand{\ben}{\begin{enumerate}}
\newcommand{\een}{\end{enumerate}}

\newcommand{\be}{\begin{equation}}
\newcommand{\ee}{\end{equation}}

\def\p{\partial}
\def\norm[#1]#2{\|#2\|_{#1}}

\def\O{B_R}
\allowdisplaybreaks
\makeatletter      
\@addtoreset{equation}{section}
\makeatother
\begin{document}

\title{Existence theorems for the Cauchy problem of 2D nonhomogeneous incompressible non-resistive MHD equations with vacuum}
\date{}
\author{ Mingtao Chen$^{\lowercase{a}}$,   Aibin Zang$^{\lowercase{b} \thanks{Corresponding author}}$}
\maketitle
\begin{center}
$^a$ School of Mathematics and Statistics, Shandong University, \\
Weihai, Weihai 264209, P.R. China.\\
Email: mtchen@sdu.edu.cn

$^b$  Center of Applied Mathematics, Yichun University, Yichun, Jiangxi, 336000, P.R.China\\
Email: 204124@jxycu.edu.cn
\end{center}

\begin{abstract}
In this paper, we investigate the Cauchy problem of the nonhomogeneous incompressible non-resistive MHD on $\mr^2$ with vacuum as far field density and prove that the 2D Cauchy problem has a unique local strong solution provided that the initial density and magnetic field decay not too slow at infinity. Furthermore, if the initial data satisfy some additional regularity and compatibility conditions, the strong solution becomes a classical one.
\end{abstract}


\textbf{Mathematics Subject Classication(2010)}: 35Q35; 74D10
\section{Introduction}
Nonhomogeneous incompressible non-resistive magnetohydrodynamics (MHD) equations  in $\mathbb{R}^2$ express as follows:
\be\la{1.1}
\left\{
\begin{array}{lll}
\n_t+(\uu\cdot\nabla)\n=0,\\
(\n \uu)_t+\div(\n \uu\otimes \uu)+\nabla \pi=\mu \triangle \uu+(\na\times \hh)\times\hh,\\
\hh_t+\uu\cdot\na\hh=\hh\cdot\na\uu,\\
\div\uu=0,\quad \div\hh=0,
\end{array}
\right.
\ee
with the initial triplet
\be\la{1.2}
(\n, \uu, \hh)(0, \xx)=(\n_0, \uu_0, \hh_0)(\xx),\quad \xx\in\mr^2,
\ee
and far field behaviors (in some weak sense)
\be\la{1.3}
\uu(t, \xx)\rightarrow \mathbf{0},\quad \n(t, \xx)\rightarrow 0,\quad \hh(t, \xx)\rightarrow 0\quad {\rm as}\ |\xx|\rightarrow\infty, \quad {\rm for}\ t\geq0.
\ee
Here $\n=\n(t,\xx)$, $\uu=(\uu_1, \uu_2)(t, \xx)$, $\hh=(\hh_1, \hh_2)(t, \xx)$ and $\pi=\pi(t, \xx)$ represent the unknown density, velocity, magnetic and pressure of the fluid, respectively. $\mu>0$ stands for the viscosity coefficient.

Magnetohydrodynamics equations describe  the motion of conducting fluids in an electromagnetic field, which has a very broad range of applications. Since the dynamic motion of the fluid and the magnetic field interact deeply on each other, the hydrodynamic and electrodynamic effects are strongly coupled, which were stated in \cite{ca, ku} and references therein.  

Note that if $\hh=0$, the MHD system \eqref{1.1} reduces to the well-known nonhomogeneous incompressible Navier-Stokes equations, which have been discussed in\cite{kazhikov, kazhikov1, cho2, choe, craig, huang, kim, liang, lv3, lions96, jianwen} and references therein, while nonhomogeneous incompressible MHD equations have been studied and obtained a few results by many mathematicians in \cite{chen, lv, lv2, sixin} and references therein. As far as the system (1.1)  is concerned, it is called the viscous and non-resistive incompressible MHD equations, which describe the conducting fluids with a very high conductivity and the third equation of the system (1.1) means that the magnetic field move along exactly with the fluid, rather that simple diffusing out. The viscous and non-resistive incompressible MHD equations are established on the  physical background in \cite{ca, ku, cha, fre} and references therein.

There are many results concerning with the multi-dimensional  incompressible MHD. For instance, as $\n\equiv$ constant, i.e., (1.1) is called homogeneous incompressible non-resistive MHD. Jiu and Niu \cite{jiu} established the local existence of solutions in 2D as the initial data in $H^s$ for integer $s\geq3$.   In the remarkable work  \cite{lin},  Lin et al. proved the existence of global-in-time solutions for initial data sufficiently close to certain equilibrium solutions in two dimensional Lagrangian coordinates by the modern analysis techniques.  Ren et al. \cite{ren} improved the results in \cite{lin} without imposing some admissible conditions for initial data and confirmed that  the energy of the MHD equations is dissipated at a rate independent of the ohmic resistivity. In 3D periodic domain, Pan et al. \cite{pan} established the the global existence of classical solutions to incompressible viscous magneto-hydrodynamical system with zero magnetic diffusion  if the initial magnetic field is close to an equilibrium state and the initial data have some symmetries. Fefferman et al. \cite{feff1}  improved that results of \cite{jiu} to d-dimensional  space and then in \cite{feff2} presented an inspiring local-in-time existence and uniqueness solutions in nearly optimal Sobolev space in $\mr^d$ ($d=2, 3$) for $\hh_0\in H^s(\mr^d)$ and $\uu_0\in H^{s-1+\ve}(\mr^d)$ with $s>d/2$ and $0<\ve<1$. It is worth noting that  Chemin et al. \cite{che} proved the local existence in Besov spaces with less regularity assumptions on $\uu_0$ than that of $\hh_0$ due to the existence of the diffusive term in the momentum equations, precisely, $\uu_0\in B_{2, 1}^{d/2-1}(\mr^d)$ and $\hh_0\in B_{2, 1}^{d/2}(\mr^d)$.

It is also an important issue to study the fluid equations with vacuum as  far field density. Recently, Li et al. \cite{lliang} proved the local existence of strong/classical solutions to the 2-D Cauchy problem of the compressible Navier-Stokes equations on the whole plane $\mr^2$ with vacuum as far field density. With the aid of \cite{lliang}, Liang \cite{liang} and Lv et al. \cite{lv, lv2, lv3}  have obtained the well-posedness for the flows with vacuum as field density, such as nonhomogeneous incompressible Navier-Stokes and  the nonhomogeneous incompressible  MHD equations with resistive term.  However, there are few results about the non-homogeneous incompressible MHD equations without resistive term, specially, the initial density may have compact support. Therefore, we try to discuss on the well-posedness of classic solutions to the system \eqref{1.1} and then obtain the following existence theorems.

%
%

\begin{theorem}\la{thm1}
Let $\eta_0$ be a positive constant and define
\be\la{1.4}
\bar\xx\triangleq (e+|\xx|^2)^{1/2}\ln^{1+\eta_0}(e+|\xx|^2),
\ee
denote $$
L^r=L^r(\mr^2), \quad W^{s, r}=W^{s, r}(\mr^2), \quad H^s=W^{s, 2},
$$ where $1\leq r\leq \infty$.
For constants $q>2$ and $a>1$, assume that the initial triplet $(\n_0, \uu_0, \hh_0)$ satisfy that
\be\la{1.5}
\left\{
\begin{array}{lll}
\n_0\geq0, \ \n_0\bar\xx^a\in L^1\cap H^1\cap W^{1, q},\ \hh_0 \bar\xx^a \in  H^1\cap W^{1, q},\\
\na\uu_0\in L^2,\ \sqrt{\n_0}\uu_0\in L^2,\ \div\uu_0=\div\hh_0=0.
\end{array}
\right.
\ee
Then there exists a positive time $T_0>0$ such that the problem \eqref{1.1}--\eqref{1.3} has a unique strong solution quaternion $(\n, \uu, \pi, \hh)$ on $\mr^2\times (0, T_0]$ with $\div\uu=0$ and $ \div\hh=0$, satisfying that
\be\la{1.6}
\left\{
\begin{array}{lll}
\n\in C([0, T_0]; L^1\cap H^1\cap W^{1, q}), \ \n\bar\xx^a\in L^\infty (0, T_0; L^1\cap H^1\cap W^{1, q}),\\
\hh\in C([0, T_0]; L^2\cap H^1\cap W^{1, q}), \ \hh\bar\xx^a\in L^\infty(0, T_0; H^1\cap W^{1, q}),\\
\sqrt\n\uu,\ \na\uu,  \ \sqrt{t}\sqrt\n\uu_t, \ \sqrt{t}\na\pi, \ \sqrt t\na^2\uu\in L^\infty (0, T_0; L^2),\\
\na\uu\in L^2(0, T_0; H^1)\cap L^{(q+1)/q}(0, T_0; W^{1, q}), \ \sqrt{t} \na\uu\in L^2(0, T_0; W^{1, q}), \\
\na\pi\in L^2(0, T_0; L^2)\cap L^{(q+1)/q}(0, T_0; L^q),\\
\sqrt\n \uu_t, \sqrt t\na\uu_t\in L^2(\mr^2\times (0, T_0)),
\end{array}
\right.
\ee
and that
\be\la{1.7}
\inf_{0\leq t\leq T_0}\int_{B_N}\n(t, \xx)d\xx\geq \frac14\int_{\mr^2} \n_0(\xx)d\xx,
\ee
for some constant $N>0$ and $B_N\triangleq\{\xx\in\mr^2||\xx|<N\}$.
\end{theorem}

Moreover, if the initial data $(\n_0, \uu_0, \hh_0)$ satisfy some additional regularity and compatibility condition, the local strong solution $(\n, \uu, \pi, \hh)$ obtained in Theorem \ref{thm1} becomes a classical one, that is,
\begin{theorem}\la{thm2}
In addition to \eqref{1.5}, suppose that
\be\la{1.8}
\na^2\n_0, \ \na^2\hh_0\in L^2\cap L^q,
\quad\bar\xx^{\delta_0}\na^2\n_0, \ \bar\xx^{\delta_0}\na^2\hh_0, \ \na^2\uu_0\in L^2,
\ee
for some constant $\delta_0\in (0, 1)$. Moreover, assume that the following compatibility condition holds for some $\mathbf{g}\in L^2$,
\be\la{1.9}
-\mu\triangle\uu_0+\na \pi_0-(\na\times \hh_0)\times\hh_0=\n_0^{1/2}\mathbf{g}.
\ee
Then, in addition to \eqref{1.6} and \eqref{1.7}, the strong solution $(\n, \uu, \pi, \hh)$ obtained in Theorem \ref{thm1} satisfies
\be\la{1.10}
\left\{
\begin{array}{lll}
\na^2\n, \ \na^2\hh\in C([0, T_0]; L^2\cap L^q),\ \bar\xx^{\delta_0}\na^2\n,  \ \bar\xx^{\delta_0}\na^2\hh \in L^\infty(0, T_0; L^2),\\
\na^2\uu, \ \sqrt\n\uu_t, \ \sqrt t\na\uu_t, \ t\sqrt\n \uu_{tt},\ t\na^2\uu_t\in L^\infty(0, T_0; L^2),\\
\na \pi,\ \sqrt t\pi_t,\ t\pi_t\in L^\infty(0, T_0; L^2),\\
t\na^3\uu, \ t\na^2\pi \in L^\infty(0, T_0; L^2\cap L^q), \\
\na\uu_t,  \ t\na\uu_{tt}\in L^2(0, T_0; L^2),\\
t\na^2(\n\uu)\in L^\infty(0, T_0; L^{(q+2)/2}).
\end{array}
\right.
\ee
\end{theorem}

As far as our existence results are concerned, it is easy to find that as $\hh=0$, Theorem \ref{thm1} is similar to the existence ressult in \cite{liang}.  While absenting the more assumptions for initial data as in \cite{lin} or \cite {ren} and excluding the higher derives of $\hh$ in the model \eqref{1.1}, it is possible not to  expect to obtain the global-in-time solution which is the same to the ones for nonhomogeneous Navier-Stokes equations \cite{lv3}, and resistive MHD equations \cite{lv}, but to establish a simple criterion depending only on $\hh$, with the motivation of the ideas in \cite{chenzang}, \cite{lv3} and \cite{zhou}, precisely,
$$
\underset{T\rightarrow T^*}{\lim\sup}\left(\|\hh\|_{L^\infty(0, T^*; L^\infty (\mr^2))}+\|\na\hh\| _{L^1(0, T_*; L^r(\mr^2)}\right)=\infty,
$$
where $r>2$ and $T^*$ is the maximal time of existence local classical solutions. For above claim, we just state it here without giving proof in detail.

In order to check the existence theorems, i.e. Theorem \ref{thm1} and Theorem \ref{thm2}, we manage to overcome the following difficulties. Firstly, in the plane, when  the far fireld density is vacuum, there is immediate difference between 2D case and 3D one, more precisely, the $L^p-$norm of $\uu$ could not be bounded directly by the $L^2-$norms of $\sqrt\n\uu$ and $\nabla\uu$. Therefore, it means that the methods  successfully used in \cite{chen, cho2, choe} can not be simply applied to our model. Utilizing the remarkable techniques in \cite{lions96, lliang, liang, lv}, we still treat the strong coupling term between the velocity vector field and the magnetic field, such as the terms $\uu\cdot\nabla\hh$ and $(\nabla\times\hh)\times\hh$ and so on. To end it,  one will borrow some ideas from \cite{lv}, and establishes the spatial weighted estimates for the magnetic fields and the density. Secondly, there is another difficulty caused by the lack of resistive term in magnetic equation $(1.1)_3$. It is impossible to infer the estimates for higher order derivatives of magnetic fields ``$\hh$" and ``$\nabla\hh$". However, due to $\div\hh=0$, let $\hh=(\partial_2\phi,-\partial_1\phi)$ for some potential $\phi.$ It is easy to see that the equation $(1.1)_3$ can be written as 
$$\phi_t+\uu\cdot\nabla\phi=0,$$
which shows that there is the same structure between the equation $(1.1)_1$ and $(1.1)_3$. Whence  the techniques to the estimates for the magnetic fields could be a homology with the one for the density, and then one can succeed in obtaining the same regular estimates for the magnetic fields and the density.

The rest of the paper is organized as follows: In Section 2, we collect some elementary facts and inequalities which will be needed in later analysis. Sections 3 and 4 are devoted to the a priori
estimates which are needed to obtain the local existence and uniqueness of strong/classical solutions. The main results: Theorem \ref{thm1} and \ref{thm2} are proved in Section 5.

\section{Preliminaries}

First, the following local existence theory on bounded balls, where the initial density is strictly away from vacuum, can be shown by similar arguments as in \cite{chen, choe, cho2}.
\bl\la{l2.1}
For $R>0$ and $B_R=\{\xx\in\mr^2||\xx|<R\}$, assume that the triplet $(\n_0, \uu_0, \hh_0)$, satisfies
\be\la{2.1}
(\n_0, \uu_0, \hh_0)\in H^3(B_R),\  \ \div\uu_0=0, \div\hh_0=0,\ \inf_{\xx\in B_R}\n_0(\xx)>0.
\ee
Then there exist a small time $T_R>0$ and a unique classical solution $(\n, \uu, \pi, \hh)$ to the following initial-boundary-value problem
\be\la{2.2}
\left\{
\begin{array}{lll}
\n_t+\div(\n\uu)=0,\\
\n\uu_t+\n(\uu\cdot\na)\uu+\na \pi=\mu\triangle\uu+(\nabla\times\hh)\times\hh,\\
\hh_t+(\uu\cdot\na)\hh=(\hh\cdot\na)\uu,\\
\div\uu=0,\quad \div\hh=0,\\
\uu =0,\xx\in\p B_R,\quad t>0,\\
(\n, \uu, \hh)(0, \xx)=(\n_0, \uu_0, \hh_0)(\xx), \quad\xx\in B_R,
\end{array}
\right.
\ee
on $B_R\times [0, T_R]$ such that
\be\la{2.3}
\left\{
\begin{array}{lll}
\n, \hh\in C([0, T_R]; H^3), \ \uu\in C([0, T_R]; H^3)\cap L^2(0, T_R; H^4), \\
\pi\in C([0, T_R]; H^2)\cap L^2(0, T_R; H^3),
\end{array}
\right.
\ee
where we denote $L^2=L^2(B_R)$ and $H^k=H^k(B_R)$ for some positive integer $k$.
\el

Then, for either $\Omega=\mr^2$ or $\Omega=\O$ with $R\geq1$, the following weighted $L^p$-bounds for elements of Hilbert space $\ti D^{1, 2} (\Omega) \triangleq \{v\in H_{\rm loc}^1(\Omega)|\na v\in L^2(\Omega)\}$ will play a crucial role in our analysis, which can be found in \cite[Lemma 2.3]{liang}.
\bl\la{l2.3}
Let $\bar{\mathbf{x}}$ be as in \eqref{1.5} and $\Omega=\mr^2$ or $\Omega=\O$ with $R\geq 1$.  Suppose that $\n\in L^1(\Omega)\cap L^\infty(\Omega)$ is a non-negative function such that
\be\la{2.4}
M_1\leq\int_{B_{N_1}}\n d\xx,\quad \|\n\|_{ L^1(\Omega)\cap L^\infty(\Omega)}\leq M_2,
\ee
for some positive constants $M_1$, $M_2$ and $N_1\geq1$ with $B_{N_1}\subset\Omega$. Then,
 for $\ve>0$ and $\eta>0$, there is a positive constant $C$ depending solely on $\ve, \eta, M_1, M_2, N_1$ and $\eta_0$ such that any $\mathbf v\in\ti D^{1, 2}(\om)$ satisfies
\be\la{2.5}
\|\mathbf{v}\bar\xx^{-\eta}\|_{L^{(2+\varepsilon)/\tilde{\eta}}(\Omega)} \leq C\|\sqrt\n\mathbf v\|_{L^2(\Omega)}+C\|\na \mathbf v\|_{L^2(\Omega)},
\ee
with $\tilde{\eta}=\min\{1,\eta\}$.
\el

Next, we consider the following Stokes system,
\be\la{2.6}
\left\{
\begin{array}{lll}
-\triangle \uu+\nabla \pi=F, &{\rm in\ }B_R,\\
\nabla\cdot u=0,&{\rm in\ }B_R\\
\uu=0, &{\rm on\ }\p B_R.
\end{array}
\right.
\ee

The proof of the following $L^p$-bound can be found in  \cite[Theorem IV.6.1]{galdi}.
\bl\la{l2.5}
Let $\uu\in W_0^{1, q}(B_R)$ be a weak solution of the system \eqref{2.6}, where $q> 1$. If $F\in W^{k, q}(B_R)$ for $k\geq 0$, then $\uu\in W^{k+2, q}(B_R)$ and
\be\la{2.7}
\|\nabla^{k+2}\uu\|_{L^{q}(B_R)}+\|\nabla^{k+1}\pi\|_{L^q}\leq C\| F\|_{W^{k, q}(B_R)},
\ee
where $C$ independent of $R$.
\el

\section{A priori estimates (I)}
In this section and the next, for $p\in[1, \infty]$ and $k\geq0$, we denote
$$
\int f d\xx=\int_{B_R} fd\xx, \quad L^p=L^p(B_R), \quad W^{k, p}=W^{k, p}(B_R), \quad H^k=W^{k, 2}.
$$
Moreover, for $R>4N_0\geq 4$, we assume that the smooth triplet $(\n_0, \uu_0, \hh_0)$ satisfies, in addition to \eqref{2.1}, that
\be\la{3.1}
\frac12\leq \int_{B_{N_0}}\n_0(\xx)d\xx\leq \int_{B_R}\n_0(\xx)d\xx\leq \frac32.
\ee
It follows from Lemma \ref{2.1} that there exists some $T_R>0$ such that the initial-boundary-value problem \eqref{2.2} has a unique classical solution $(\n, \uu, \pi, \hh)$ on $B_R\times [0, T_R]$ satisfying \eqref{2.3}. For $\bar\xx, \eta_0, a$ and $q$ as in Theorem \ref{thm1}, the main aim of this section is to derive the following key a priori estimate on $\phi(t)$ defined by
\be\la{3.2}
\phi(t)\triangleq 1+\|\sqrt\n\uu\|_{L^2}+\|\na\uu\|_{L^2}+\|\hh\|_{L^2} +\|\hh\bar\xx^a\|_{H^1\cap W^{1, q}}+\|\n\bar\xx^a\|_{L^1\cap H^1\cap W^{1, q}}.
\ee

\begin{proposition}\la{prop}
Assume that $(\n_0, \uu_0, \hh_0)$ satisfies \eqref{2.1} and \eqref{3.1}. Let $(\n, \uu, \pi, \hh)$ be the solution to the initial-boundary-value problem \eqref{2.2} on $B_R\times (0, T_R]$ obtained by Lemma \ref{2.1}. Then there exist positive constants $T_0$ and $M$ both depending solely on $\mu, q, a, \eta_0, N_0$ and $C_0$ such that
\be\la{3.3}
\sup_{0\leq t\leq T_0}\phi(t)+\int_0^{T_0}\left(\|\na^2 \uu\|_{L^q} ^{(q+1)/q}+t\|\na^2 \uu\|_{L^q}^2+\|\na^2 \uu\|_{L^2}^2\right)dt\leq M,
\ee
where
$$
C_0=\|\sqrt{\n_0}\uu_0\|_{L^2}+\|\na\uu_0\|_{L^2}+\|\hh_0\|_{L^2} +\|\hh_0\bar\xx^a\|_{H^1\cap W^{1, q}}+\|\n_0\bar\xx^a\| _{L^1\cap H^1\cap W^{1, q}}.
$$
\end{proposition}

The proof of Proposition \ref{prop} will be postponed at the end of this section. First, we start with the following energy estimate for $(\n, \uu, \pi, \hh)$ and preliminary $L^2$-bounds for $\na\uu$.
\bl\la{l3.1}
Let $(\n, \uu, \pi, \hh)$ be a smooth solution to the initial-boundary-value problem \eqref{2.2}. Then there exist a positive constant $\alpha=\alpha( q)>1$ and a $T_1=T_1(C_0, N_0)>0$ such that for all $t\in (0, T_1]$,
\begin{align}\la{3.4}
&\sup_{0\leq s\leq t}\left(\|\na\uu\|_{L^2}^2+\|\hh\|_{L^2}^2+ \|\sqrt\n\uu\|_{L^2}^2 \right)\\
&\qquad\quad+\int_0^t\left(\|\na\uu\|_{L^2}^2+\|\sqrt\n \uu_t\|_{L^2}^2\right)ds
\leq C+C\int_0^t\phi^\alpha(s) ds.\nonumber
\end{align}
\el
\pf Firstly, by energy estimates we obtain
\be\la{3.5}
\sup_{0\leq s\leq t}\left(\|\sqrt\n\uu\|_{L^2}^2+\|\hh\|_{L^2}^2 \right) +\int_0^t\|\na\uu\|_{L^2}^2ds\leq C,\qquad
\sup_{0\le s\le t}\|\n\|_{L^1\cap L^\infty}\le C.
\ee

Next, for $N>1$ and $\vp_N\in C_0^\infty(B_R)$ such that
\be\la{3.6}
0\leq\vp_N\leq 1, \ \vp_N(\xx)=1,\ {\rm if}\ |\xx|\leq N/2,\ \ |\na^k\vp_N|\leq CN^{-k} (k=1, 2),
\ee
then it follows from \eqref{3.1} and \eqref{3.5} that
\begin{align}\la{3.7}
\frac{d}{dt}\int\n\vp_{2N_0}d\xx=&\int\n\uu\cdot\na\vp_{2N_0}d\xx\\
\geq&-CN_0^{-1}\left(\int\n d\xx\right)^{1/2}\left( \int\n|\uu|^2 \right)^{1/2}\geq -\ti C(C_0, N_0),\nonumber
\end{align}
where in the last inequality we have used
$$
\int\n d\xx=\int\n_0 d\xx,
$$
due to \eqref{2.2}$_1$ and \eqref{2.2}$_4$. Integrating \eqref{3.7} over $(0, T_1)$ shows
\be\la{3.8}
\inf_{0\leq t\leq T_1}\int_{B_{2N_0}}\n d\xx\geq \inf_{0\leq t\leq T_1} \int\n\vp_{2N_0}d\xx\geq \int\n_0\vp_{2N_0}d\xx-\ti CT_1\geq1/4,
\ee
where $T_1\triangleq\min\{1, (4\ti C)^{-1}\}$. From now on, we will always suppose that $t\leq T_1$. The combination of \eqref{2.5}, \eqref{3.5} and \eqref{3.8} shows that for $\ve>0$ and $\eta>0$, every $\mathbf v\in \ti D^{1, 2}(B_R)$ satisfies
\be\la{3.9}
\|\mathbf v\bar\xx^{-\eta}\|_{L^{(2+\ve)/\ti\eta}}^2\leq C(\ve, \eta)\|\sqrt\n\mathbf v\|_{L^2}^2+C(\ve, \eta)\|\na\mathbf v\|_{L^2}^2,
\ee
with $\ti\eta=\min\{1, \eta\}$.

Multiplying \eqref{2.2}$_1$ by $\bar\xx^a$  and integrate by parts, from \eqref{3.9} for $\mathbf u\in \ti D^{1, 2}(B_R)$, note that
\begin{align*}
\frac{d}{dt}\int\n\bar\xx^ad\xx\leq &C\int \n |\uu|\bar\xx^{a-1} \ln^{1+\eta_0}(e+|\xx|^2)d\xx\\
\leq &C\|\n\bar\xx^{a-1+8/(8+a)}\|_{L^{(8+a)/(7+a)}}\|\uu \bar\xx^{-4/(8+a)}\|_{L^{8+a}}\\
\leq &C\phi^\alpha(t),
\end{align*}
which implies
\be\la{3.10a}
\sup_{0\leq s\leq t}\|\n\bar\xx^a\|_{L^1}\leq C\exp\left(C \int_0^t\phi^\alpha(s)ds\right).
\ee
By calculation from \eqref{3.10a} and \eqref{3.9}, we conclude that
\be\la{3.10}
\|\n^\eta\uu\|_{L^{(2+\ve)/\ti\eta}}+\|\uu\bar\xx^{-\eta}\| _{L^{(2+\ve)/\ti\eta}}\leq C(\ve, \eta)\phi^{1+\eta}(t).
\ee

Next, multiplying \eqref{2.2}$_2$ by $\uu_t$ and integration by parts, one yields
\begin{align}\la{3.11}
\frac{d}{dt}\left(\mu\|\na\uu\|_{L^2}^2\right)+\|\sqrt\n\uu_t\|_{L^2}^2
\leq C\int \n|\uu|^2|\na\uu|^2d\xx+2\int (\hh\cdot\na)\hh \cdot\uu_t d\xx.
\end{align}
Now we need to estimate each term on the right-hand side of \eqref{3.11}. At first, the Gagliardo-Nirenberg inequality implies that for all $p\in(2, +\infty)$,
\be\la{3.12}
\|\na\uu\|_{L^p}\leq C\|\na\uu\|_{L^2}^{2/p}\|\na\uu\|_{H^1}^{1-2/p} \leq C\phi(t)+C\phi(t)\|\na^2\uu\|_{L^2}^{1-2/p},
\ee
which together with \eqref{3.10}, it follows that for $\eta>0$ and $\ti\eta=\min\{1, \eta\}$,

\begin{align}\la{3.13}
\int\n^\eta|\uu|^2|\na\uu|^2d\xx\leq & C\|\n^{\eta/2}\uu\|_{L^{8/\ti\eta} }^2\|\na\uu\|_{L^{8/(4-\ti\eta)}}^2\\
\leq& C(\eta)\phi^{4+2\eta}(t)\left(1+\|\na^2\uu\|_{L^2}^{\ti\eta/2} \right)\no
\leq&C\phi^{\alpha(\eta)}(t)+\ve\|\na^2\uu\|_{L^2}^2.\nonumber
\end{align}

Then,  integrating by parts  with  respect to the variable $t$ and from \eqref{2.2}$_3$, one obtains
\begin{align}\la{3.16}
2\int (\hh\cdot\na)\hh \cdot\uu_t d\xx=&-2\frac{d}{dt}\int (\hh\cdot\na ) \uu\cdot\hh d\xx+2\int (\hh_t\cdot\na)\uu\cdot\hh d\xx\\
&\quad+2\int (\hh\cdot\na )\uu\cdot\hh_t d\xx\no
=&-2\frac{d}{dt}\int (\hh\cdot\na ) \uu\cdot\hh d\xx+2\int ((\hh\cdot\na)\uu\cdot\na)\uu\cdot\hh d\xx\no
&-2\int ((\uu\cdot\na)\hh\cdot\na)\uu\cdot\hh d\xx\no
&+2\int (\hh\cdot\na )\uu\cdot(\hh\cdot\na)\uu d\xx-2\int (\hh\cdot\na )\uu\cdot(\uu\cdot\na)\hh d\xx.\nonumber
\end{align}
First, it is easy to see that
\begin{align*}
\left|\int ((\hh\cdot\na)\uu\cdot\na)\uu\cdot\hh d\xx\right|+\left|\int (\hh\cdot\na )\uu\cdot(\hh\cdot\na)\uu d\xx\right|\leq C\|\hh\|_{L^\infty}^2\|\na\uu\|_{L^2}^2.
\end{align*}
Next, by H\"{o}lder inequality, Young inequality, and \eqref{3.12}, one observes
\begin{align*}
\left|\int ((\uu\cdot\na)\hh\cdot\na)\uu\cdot\hh d\xx\right| &+\left|\int (\hh\cdot\na )\uu\cdot(\uu\cdot\na)\hh d\xx\right|\no
\leq &C\int|\hh||\na\uu||\uu||\na\hh|d\xx\\
\leq &C\|\hh\|_{L^\infty}\|\uu\bar\xx^{-1}\|_{L^{2a}}\|\hh\bar\xx^a\| _{H^1}\|\na\uu\|_{L^{2a/(a-1)}}\no
\leq &C\phi^\alpha(t)+\ve\phi^{-1}(t)\|\na^2\uu\|_{L^2}^2.
\end{align*}
Substituting the above two estimates into \eqref{3.16}, we have given
\be\la{3.17}
2\int (\hh\cdot\na)\hh \cdot\uu_t d\xx\leq-2\frac{d}{dt}\int (\hh\cdot\na ) \uu\cdot\hh d\xx +C\phi^\alpha(t)+\ve\phi^{-1}(t)\|\na^2\uu\|_{L^2}^2.
\ee

Inserting \eqref{3.13} and \eqref{3.17} into \eqref{3.11}, one shows
\begin{align}\la{3.19}
&\frac{d}{dt}\mu\|\na\uu\|_{L^2}^2
+2\frac{d}{dt}\int \left( (\hh\cdot\na ) \uu\cdot\hh\right) d\xx+\|\sqrt\n\uu_t\|_{L^2}^2\\
\leq &C\phi^\alpha(t)+4\ve\phi^{-1}(t)\|\na^2\uu\|_{L^2}^2.\nonumber
\end{align}

To estimate the last term on the right-hand side of \eqref{3.19}, it follows from \eqref{2.7} with $F=\n\uu_t+\n\uu\cdot\na\uu+\hh\cdot\na \hh$ that for $p\in [2, q]$,
\begin{align}\la{3.20}
\|\na^2\uu\|_{L^p}\leq & C\left(\|\n\uu_t\|_{L^p}+\|\n\uu\cdot\na\uu\|_{L^p}+\||\hh||\na \hh|\|_{L^p}\right).
\end{align}
 Combining with \eqref{3.20}, \eqref{3.12} and \eqref{3.13}, we yield
\begin{align}\la{3.21}
\|\na^2\uu\|_{L^2}\leq& C\phi^{1/2}(t)\|\sqrt\n\uu_t\|_{L^2} +C\|\n\uu\cdot\na\uu\|_{L^2} +C\phi^\alpha (t)\\
\leq &C\phi^{1/2}(t)\|\sqrt\n\uu_t\|_{L^2}+\frac12\|\na^2\uu\|_{L^2} +C\phi^\alpha (t).\nonumber
\end{align}
Substituting \eqref{3.21} into \eqref{3.19}, then integrating the resultant inequality over $(0, t)$, and choosing $\ve$ suitably small, let us lead to
\begin{align}\la{3.22}
\frac\mu2\|\na\uu\|_{L^2}^2+\int_0^t \|\sqrt\n\uu_t\|_{L^2}^2ds
&\leq C+C\|\hh\|_{L^4}^4+C\int_0^t\phi^\alpha(s)ds,
\end{align}
where in the last inequality we have used
\begin{equation*}
\int \left( (\hh\cdot\na ) \uu\cdot\hh\right) d\xx\ge-\frac{\mu}{2}\|\na\uu\|_{L^2}^2-C\|\hh\|_{L^4}^4.
\end{equation*}

To estimate the second term on the right-hand side of \eqref{3.22}, multiplying \eqref{2.2}$_3$ by $4|\hh|^2\hh$ and integrating the resultant equality over $B_R$, we have
\begin{align*}
\frac{d}{dt}\|\hh\|_{L^4}^4\leq& C\int|\na\uu||\hh|^4d\xx\\
\leq &C\|\na\uu\|_{L^2}\|\hh\|_{L^2}\|\hh\|_{L^\infty}^3.\nonumber
\end{align*}
Integrating the above inequality over $(0, t)$, one infers
\be\la{3.23}
\|\hh\|_{L^4}^4\leq C+\int_0^t\phi^\alpha(s)ds.
\ee
Putting \eqref{3.23} into \eqref{3.22}, together with \eqref{3.5}, one leads to \eqref{3.4}. Therefore, we complete the proof of Lemma \ref{l3.1}. \hfill $\Box$

\bl\la{l3.2}
Let $(\n, \uu, \pi, \hh)$ and $T_1$ be as in Lemma \ref{l3.1}. Then for all $t\in (0, T_1]$,
\be\la{3.24}
\sup_{0\leq s\leq t}s\|\sqrt\n\uu_t\|_{L^2}^2+\int_0^ts\|\na\uu_t\| _{L^2}^2ds\leq C\exp\left(C\int_0^t\phi^\alpha(s)ds\right).
\ee
\el
\pf Differentiating \eqref{2.2}$_2$ with respect to $t$ yields
\begin{align}\la{3.25}
&\n\uu_{tt}+\n\uu\cdot\na\uu_t-\mu\triangle\uu_t \\
=&-\n_t(\uu_t+\uu\cdot\na\uu)-\n\uu_t\cdot\na\uu-\na \pi_t+\hh_t\cdot\na\hh+\hh\cdot\na\hh_t.\nonumber
\end{align}
Multiplying \eqref{3.25} by $\uu_t$ and integrating the resultant equation over $B_R$, we obtain
\begin{align}\la{3.26}
&\frac12\frac{d}{dt}\int\n|\uu_t|^2d\xx+\mu\int|\na\uu_t|^2d\xx\\
=&-2\int\n\uu\cdot\na\uu_t\cdot\uu_t d\xx-\int\n\uu\cdot\na(\uu\cdot \na\uu\cdot\uu_t)d\xx\no
&-\int\n\uu_t\cdot\na\uu\cdot\uu_t d\xx-\int(\hh\otimes\hh)_t:\na\uu_td\xx\no
\leq &C\int\n|\uu||\uu_t|\left(|\na\uu_t|+|\na\uu|^2+|\uu||\na^2\uu |\right)d\xx+C\int\n|\uu|^2|\na\uu||\na\uu_t|d\xx\no
&+C\int\n|\uu_t|^2|\na\uu|d\xx+C \int |\hh||\hh_t||\na\uu_t|d\xx.\nonumber
\end{align}

We now estimate each term on the right-hand side of \eqref{3.26} as follows:

First, it follows from \eqref{3.2}, \eqref{3.5}, \eqref{3.9}, \eqref{3.10} and \eqref{3.12} that for $\ve\in (0, 1)$,
\begin{align}\la{3.27}
&\int\n|\uu||\uu_t|\left(|\na\uu_t|+|\na\uu|^2+|\uu||\na^2\uu |\right)d\xx\\
\leq &C\|\sqrt\n\uu\|_{L^6}\|\sqrt\n\uu_t\|_{L^2}^{1/2} \|\sqrt\n\uu_t\|_{L^6}^{1/2}\left(\|\na\uu_t\|_{L^2} +\|\na\uu\|_{L^4}^2\right)\no
&+C\|\n^{1/4}\uu\|_{L^{12}}^2\|\sqrt\n\uu_t\|_{L^2}^{1/2}  \|\sqrt\n\uu_t\|_{L^6}^{1/2}\|\na^2\uu\|_{L^2}\no
\leq &C\phi^\alpha(t)\|\sqrt\n\uu_t\|_{L^2}^{1/2}\left( \|\sqrt\n\uu_t\|_{L^2}^{1/2}+ \|\na\uu_t\|_{L^2}^{1/2}\right) \left(\|\na\uu_t\|_{L^2}+\|\na^2\uu\|_{L^2}+\phi(t)\right)\no
\leq&\ve\|\na\uu_t\|_{L^2}^2+C\phi^\alpha(t)\left( \|\na^2\uu\|_{L^2}^2+\|\sqrt\n\uu_t\|_{L^2}^2+1\right).\nonumber
\end{align}

From  H\"{o}lder inequality, \eqref{3.10} and \eqref{3.12}, it is easy to show that
\begin{align}\la{3.28}
\int\n|\uu|^2|\na\uu||\na\uu_t|d\xx\leq & C\|\sqrt\n\uu\|_{L^8}^2\|\na \uu\|_{L^4}\|\na\uu_t\|_{L^2}\\
\leq &\ve\|\na\uu_t\|_{L^2}^2+C\left(\phi^\alpha (t)+\|\na^2\uu\|_{L^2}^2 \right).\nonumber
\end{align}

Then, by \eqref{3.9} and H\"{o}lder inequality, we have the following estimate
\begin{align}\la{3.29}
\int\n|\uu_t|^2|\na\uu|d\xx\leq &\|\na\uu\|_{L^2}\|\sqrt\n \uu_t\|_{L^6}^{3/2} \|\sqrt\n\uu_t\|_{L^2}^{1/2}\\
\leq &\ve\|\na\uu_t\|_{L^2}^2+C\phi^\alpha(t)\|\sqrt\n \uu_t\|_{L^2}^2.\nonumber
\end{align}

Finally, due to \eqref{2.2}$_3$ and \eqref{3.10}, one obtains that
\begin{align}\la{3.31}
\int |\hh||\hh_t||\na\uu_t|d\xx
\leq &C\int \left(|\hh||\na\uu|+|\na\hh||\uu|\right)|\na\uu_t|d\xx\\
\leq &C\left(\|\hh\|_{L^\infty}\|\na\uu\|_{L^2}+\|\hh\bar\xx^a\|_{W^{1, q}} \|\uu\bar\xx^{-a}\|_{L^{2q/(q-2)}}\right)\|\na\uu_t\|_{L^2}\no
\leq &\ve\|\na\uu_t\|_{L^2}^2+C\phi^\alpha(t).\nonumber
\end{align}

Inserting \eqref{3.27}--\eqref{3.31} into \eqref{3.26} and choosing $\ve$ suitably small, we observe that
\begin{align}\la{3.32}
&\frac{d}{dt}\int\n|\uu_t|^2d\xx+\mu\int|\na\uu_t|^2d\xx\\
\leq &C\phi^\alpha(t)\left(1+\|\sqrt\n\uu_t\|_{L^2}^2 +\|\na^2\uu\|_{L^2}^2 \right)\no
\leq &C\phi^\alpha(t)\|\sqrt\n\uu_t\|_{L^2}^2+C\phi^\alpha(t),\nonumber
\end{align}
where in the last inequality we have used \eqref{3.21}. Then, multiplying \eqref{3.32} by $t$, we finally obtain \eqref{3.24} after using Gronwall inequality and \eqref{3.4}. Therefore, we complete the proof of Lemma \ref{l3.2}. \hfill $\Box$

\bl\la{l3.3}
Let $(\n, \uu, \pi, \hh)$ and $T_1$ be as in Lemma \ref{l3.1}. Then for all $t\in (0, T_1]$,
\be\la{3.33}
\sup_{0\leq s\leq t}\left(\|\n\bar\xx^a\|_{L^1\cap H^1\cap W^{1, q}}+\|\hh\bar\xx^a\|_{H^1\cap W^{1, q}}\right) \leq \exp\left(C\exp\left(C\int_0^t\phi^\alpha(s)ds\right)\right).
\ee
\el
\pf As we known, the estimate for $\sup_{0\leq s\leq t}\|\n\bar\xx^a\|_{L^1}$ was given by \eqref{3.10a} and the other estimates for $\n\bar\xx^a$ are similar to the one for $\hh\bar\xx^a$.

Next, it follows from the Sobolev inequality and \eqref{3.10} that for $0<\delta<1$,
\begin{align}\la{3.35}
\|\uu\bar\xx^{-\delta}\|_{L^\infty}\leq &C\left(\|\uu\bar\xx^{-\delta} \|_{L^{4/\delta}}+\|\na(\uu\bar\xx^{-\delta})\|_{L^3} \right)\\
\leq& C\left(\|\uu\bar\xx^{-\delta} \|_{L^{4/\delta}}+ \|\na\uu\|_{L^3}+\|\uu\bar\xx^{-\delta} \|_{L^{4/\delta}}\|\bar\xx^{-1}\na\bar\xx\|_{L^{12/(4-3\delta)}} \right)\no
\leq &C\left(\phi^\alpha (t) +\|\na^2\uu\|_{L^2}\right).\nonumber
\end{align}

Then, one derives from \eqref{2.2}$_3$ that $\ti\hh\triangleq \hh\bar\xx^a$ satisfies
\be\la{3.36}
\ti\hh_t+\uu\cdot\na \ti\hh-a\uu\cdot\na\ln\bar\xx\cdot\ti\hh =\ti\hh\cdot\na\uu.
\ee
Applying \eqref{3.35} and Gagliardo-Nirenberg inequality to the equation above,  it is easy to show that
\begin{align}\la{3.37}
\frac{d}{dt}\|\ti\hh\|_{L^2}\leq &C\left(\|\na\uu\|_{L^\infty} +\|\uu\cdot\na\ln\bar\xx\|_{L^\infty}\right)\|\ti\hh\|_{L^2}\\
\leq &C\left(\phi^\alpha(t)+\|\na^2\uu\|_{L^2\cap W^{1, q}}\right)\|\ti\hh\|_{L^2}.\nonumber
\end{align}
Moreover,  we can  infer that for $p\in [2, q]$
\begin{align}\la{3.38}
&\frac{d}{dt}\|\na\ti\hh\|_{L^p}\\
\leq & C\left( 1+\|\na\uu\|_{L^\infty} +\|\uu\cdot\na\ln\bar\xx\|_{L^\infty}\right)\|\na\ti\hh\|_{L^p}\no
&+C\left(\||\na\uu||\na\ln\bar\xx|\|_{L^p}+\||\uu||\na^2\ln\bar\xx|\| _{L^p} +\|\na^2\uu\|_{L^p}\right)\|\ti\hh\|_{L^\infty}\no
\leq &C\left(\phi^\alpha(t)+\|\na\uu\|_{L^2\cap W^{1, q}}\right)\|\na\ti\hh\|_{L^p}\no
&+C\left(\|\na\uu\|_{L^p}+\|\uu\bar\xx^{-1/4}\|_{L^\infty} \|\bar \xx^{-3/2}\|_{L^p}+\|\na^2\uu\|_{L^p}\right)\|\ti\hh\|_{L^\infty}\no
\leq& C\left(\phi^\alpha(t)+\|\na^2\uu\|_{L^2\cap L^p}\right) \left(1+\|\na\ti\hh\|_{L^p}+\|\na\ti\hh\|_{L^q}\right).\nonumber
\end{align}
Combining \eqref{3.37} and \eqref{3.38} yields
\begin{align}\la{3.39}
&\frac{d}{dt}\left(\|\ti\hh\|_{L^2}+\|\na\ti\hh\|_{L^p}\right)\\
\leq & C\left(\phi^\alpha(t)+\|\na^2\uu\|_{L^2}+\|\na^2\uu\|_{ L^p}\right) \left(1+\|\na\ti\hh\|_{L^p}+\|\ti\hh\|_{L^2}+ \|\na\ti\hh\|_{L^q}\right).\nonumber
\end{align}

As in formula (3.55) and (3.57) in \cite{lv}, we claim that
\be\la{3.41}
\begin{aligned}
\int_0^t\left(\|\na^2\uu\|_{ L^q}^{(q+1)/q}+\|\na\pi\|_{ L^q}^{(q+1)/q}\right)ds\le C\exp\left(C\int_0^t\phi^\alpha(s)ds\right),\\
\int_0^t\left(s\|\na^2\uu\|_{L^2\cap L^q}^2+s\|\na\pi\|_{L^2\cap L^q}^2\right)ds\leq C\exp\left(C\int_0^t\phi^\alpha(s)ds\right).
\end{aligned}
\ee
Substituting \eqref{3.41} into \eqref{3.39},  and  then Gronwall inequality implies \eqref{3.33}. Therefore, we complete the proof of Lemma \ref{l3.3}. \hfill $\Box$

Now, we complete the proof of Proposition \ref{prop}, which is a direct consequence of Lemma \ref{l3.1}, \ref{l3.2} and \ref{l3.3}.

{\it Proof of Proposition \ref{prop}.} It follows from \eqref{3.4} and \eqref{3.33} that
$$
\phi(t)\leq \exp\left(C\exp\left(C\int_0^t \phi^\alpha(s)ds\right)\right).
$$
Standard arguments thus show that for $\ti M\triangleq e^{Ce}$ and $T_0 \triangleq \min\{T_1, (CM^\alpha)^{-1}\}$,
$$
\sup_{0\leq t\leq T_0}\phi(t)\leq \ti M,
$$
which together with \eqref{3.4}, \eqref{3.21} and \eqref{3.41} leads to \eqref{3.3}. Then the proof of Proposition \ref{prop} is finished. \hfill $\Box$

\section{A priori estimates (II)}
In this section, in addition to $\mu,  q, a, \eta_0, N_0$ and $C_0$, the generic positive constant $C$ may depend on $\delta_0$, $\|\na^2\uu_0\|_{L^2}$, $\|\na^2\n_0\|_{L^q}$, $\|\na^2\hh_0\|_{L^q}$, $\|\bar\xx^{\delta_0}\na^2\n_0\|_{L^2}$,  $\|\bar\xx^{\delta_0}\na^2\hh_0\|_{L^2}$ and $\|\mathbf{g}\|_{L^2}$.

\bl\la{l4.1}
It holds that
\be\la{4.1}
\sup_{0\leq t\leq T_0}\left(\|\bar\xx^{\delta_0}\na^2\n\|_{L^2} +\|\bar\xx^{\delta_0} \na^2\hh\|_{L^2}\right)\leq C.
\ee
\el
\pf First, due to \eqref{1.9}, \eqref{2.1} and \eqref{2.2}$_2$, defining $$
\sqrt\n\uu_t(t=0, \xx)\triangleq -\mathbf{g}-\sqrt{\n_0}\uu_0\cdot\na\uu_0,
$$
integrating \eqref{3.32} over $(0, T_0)$ and using \eqref{3.3} and \eqref{3.4}, we have
\be\la{4.2}
\sup_{0\leq t\leq T_0}\|\sqrt\n\uu_t\|_{L^2}^2+\int_0^{T_0} \|\na\uu_t\|_{L^2}^2dt\leq C,
\ee
which together with \eqref{3.3} and \eqref{3.21} leads to
\be\la{4.3}
\sup_{0\leq t\leq T_0}\|\na\uu\|_{H^1}\leq C.
\ee
Combined with \eqref{3.3}, \eqref{3.35} and \eqref{4.3}, it is sure to show that for $\delta\in (0, 1]$,
\be\la{4.4}
\|\n^\delta\uu\|_{L^\infty}+\|\bar\xx^{-\delta}\uu\|_{L^\infty}\leq C(\delta).
\ee
Following the direct calculations, one yields that for $2\leq r\leq q$
\be\la{4.5}
\|\n_t(\bar\xx^{(1+a)/2}+|\uu|)\|_{L^r}+\|\hh_t(\bar\xx^{(1+a)/2} +|\uu|)\|_{L^r}\leq C.
\ee
Due to \eqref{2.2}$_1$, \eqref{2.2}$_3$, \eqref{3.3}, \eqref{4.3} and \eqref{4.4}, it follows from \eqref{3.9}, \eqref{4.2}, \eqref{4.3} and \eqref{4.4} that for $\delta\in (0, 1]$ and $s>2/\delta$,
\begin{align}\la{4.6}
\|\bar\xx^{-\delta}\uu_t\|_{L^s}+\|\bar\xx^{-\delta}\uu\cdot \na\uu\|_{L^s}\leq & C\|\bar\xx^{-\delta}\uu_t\|_{L^s} +C\|\bar\xx^{-\delta}\uu\|_{L^\infty}\|\na\uu\|_{L^s}\\
\leq &C(\delta, s)+C(\delta, s)\|\na\uu_t\|_{L^2}.\nonumber
\end{align}

Next, denoting  $\bar\hh\triangleq \bar\xx^{\delta_0}\hh$ and $\bar\n\triangleq \bar\xx^{\delta_0}\n$, we easily get from \eqref{3.3} that
\be\la{4.7}
\|\bar\n\|_{L^\infty}+\|\na\bar \n\|_{L^2\cap L^q}+\|\bar\hh\|_{L^\infty}+ \|\na \bar\hh\|_{L^2\cap L^q}\leq C,
\ee
where $\bar\n$ and $\bar\hh$ satisfy,
\be\la{a4.7}
\bar\n_t+\uu\cdot\na\bar\n-\delta_0 \bar\n\uu\cdot\na \ln\bar\xx=0,
\ee
and
\be\la{4.8}
\bar\hh_t+\uu\cdot\na\bar\hh-\delta_0\uu\cdot\na\ln\bar\xx\cdot \bar\hh=\bar\hh\cdot\na\uu,~\mbox{respectively.}
\ee
Therefore, by energy method, we should give that
\begin{align}\la{4.9}
\frac{d}{dt}\|\na^2\bar\hh\|_{L^2}\leq & C\left(1+\|\na\uu\| _{L^\infty} +\|\uu\cdot\na\ln\bar\xx\|_{L^\infty}\right)\|\na^2\bar\hh\|_{L^2} +C\||\na^2\uu||\na\bar\hh|\|_{L^2}\\
&+C\||\na\bar\hh||\na\uu||\na\ln\bar\xx|\|_{L^2}+C\||\na\bar\hh||\uu| |\na^2 \ln\bar\xx|\|_{L^2}\no
&+C\|\bar\hh\|_{L^\infty}\left(\|\na^2\left(\uu\cdot\na\ln\bar\xx \right)\|_{L^2}+\|\na^3\uu\|_{L^2}\right)\no
\leq &C\left(1+\|\na\uu\|_{L^\infty}\right)\|\na^2\bar\hh\|_{L^2} +C\|\na^2 \uu\|_{L^{2q/(q-2)}}\|\na\bar\hh\|_{L^q} \no
&+C\|\na\bar\hh\|_{L^2}\|\na\uu\|_{L^\infty}+C\|\na\bar\hh\|_{L^2} \||\uu||\na^2\ln\bar\xx|\|_{L^\infty}\no
&+C\|\na^2\uu\|_{L^2}+C\|\na\uu\|_{L^2}+C\||\uu||\na^3\ln\bar\xx|\|_{L^2} +C\|\na^3\uu\|_{L^2}\no
\leq & C\left(1+\|\na\uu\|_{L^\infty}\right)\|\na^2\bar\hh\|_{L^2}+C +C\|\na^3\uu\|_{L^2},\nonumber
\end{align}
here the formula \eqref{4.4} and \eqref{4.7} were used in the second and third inequalities. Similarly, we can also obtain from \eqref{a4.7}
\be\la{4.10}
\frac{d}{dt}\|\na^2\bar\n\|_{L^2}\leq C\left(1+\|\na\uu\|_{L^\infty}\right)\|\na^2\bar\n\|_{L^2}+C +C\|\na^3\uu\|_{L^2}.
\ee
Combing \eqref{4.9} with \eqref{4.10}, we get
\begin{align}\la{4.12}
&\frac{d}{dt}\left(\|\na^2\bar\n\|_{L^2}+\|\na^2\bar\hh\|_{L^2}\right) \\
\leq& C\left(1+\|\na\uu\|_{L^\infty}\right)\left(\|\na^2\bar\n\|_{L^2}+ \|\na^2\bar\hh\|_{L^2}\right) +C+C\|\na^3\uu\|_{L^2}.\nonumber
\end{align}
By  \eqref{2.7} and \eqref{4.3}, one shows that
\begin{align}\la{4.13}
\|\na^3\uu\|_{L^2}\leq &C \|\na(\n\uu_t)\|_{L^2}+C\|\na(\n\uu\cdot \na\uu)\|_{L^2}\\
&+C\|\na^2|\hh|^2\|_{L^2}+C \|\na(\hh\cdot\na\hh)\|_{L^2}\no
\leq &C\|\bar\n\|_{L^\infty}\|\na\uu_t\|_{L^2}+C\|\bar\xx^a\na\n\|_{L^q}\|\bar \xx^{-a}\uu_t\|_{L^{2q/(q-2)}}\no
&+C\|\bar\xx^{-\delta}\uu\|_{L^\infty} \|\bar\xx^\delta\na\n\|_{L^q} \|\na\uu\|_{L^{2q/(q-2)}}+C\|\bar\n\|_{L^\infty}\|\na\uu\|_{L^4}^2\no
&+C\|\bar\n\|_{L^\infty}\|\bar\xx^{-\delta}\uu\|_{L^\infty}\|\na^2\uu\| _{L^2}^2 \no
&+C\|\na\bar\hh\|_{L^4}^2 +C\|\bar\hh\|_{L^\infty}\|\na^2
\hh\|_{L^2}\no
\leq &C\|\na\uu_t\|_{L^2}+C\left(\|\na^2\bar\n\|_{L^2} +\|\na^2\bar\hh\|_{L^2}\right)+C,\nonumber
\end{align}
where in the last inequality we have used \eqref{3.3}, \eqref{4.3}, \eqref{4.4}, \eqref{4.6} and the following fact:
\begin{align}\la{4.14}
\|\bar\xx^{\delta_0}\na^2\n\|_{L^2}+\|\bar\xx^{\delta_0}\na^2 \hh\|_{L^2}
\leq C\|\na^2(\bar\xx^{\delta_0}\n)\|_{L^2}+\|\na^2(\bar\xx^{\delta_0} \hh)\|_{L^2}+C.
\end{align}
Substituting \eqref{4.13} into \eqref{4.12} , one has
\begin{align*}
&\frac{d}{dt}\left(\|\na^2(\bar\xx^{\delta_0}\n)\|_{L^2} +\|\na^2(\bar\xx^{\delta_0}\hh)\|_{L^2}\right)\no
\leq &C(1+\|\na^2\uu\|_{L^q})\left(\|\na^2(\bar\xx^{\delta_0}\n)\|_{L^2} +\|\na^2(\bar\xx^{\delta_0}\hh)\|_{L^2}\right)\no
&+C\|\na\uu_t\|_{L^2}+C,
\end{align*}
then by \eqref{3.3}, \eqref{4.2}, \eqref{4.14} and Gronwall inequality, we have got \eqref{4.1} and completes the proof of Lemma \ref{l4.1}. \hfill $\Box$

\bl\la{l4.2}
It holds that
\be\la{4.15}
\sup_{0\leq t\leq T_0}t\|\na\uu_t\|_{L^2}^2+\int_0^{T_0}t\left( \|\sqrt\n \uu_{tt}\|_{L^2}^2+\|\na^2\uu_{t}\|_{L^2}^2\right)dt\leq C.
\ee
\el
\pf Multiplying \eqref{3.25} by $\uu_{tt}$ and integrating the resultant equality over $B_R$,  integration by parts then we lead to
\begin{align}\la{4.16}
&\frac12\frac{d}{dt}\left(\mu\|\na\uu_t\|_{L^2}^2\right)+\|\sqrt\n \uu_{tt}\|_{L^2}^2\\
=&-\int\left(2\n\uu\cdot\na\uu_t\cdot\uu_{tt}+\n\uu_t\cdot\na \uu \cdot\uu_{tt}\right)d\xx-\int\n\uu\cdot\na(\uu\cdot\na\uu)\cdot \uu_{tt}d\xx\no
&-\int\n\uu\cdot\na\uu_{tt}\cdot\uu_td\xx-\int\n\uu\cdot\na\uu_{tt} \cdot(\uu \cdot\na)\uu d\xx\no
&-\int\hh_t\cdot\na\uu_{tt}\cdot \hh d\xx-\int\hh\cdot\na\uu_{tt}\cdot\hh_td\xx.\nonumber
\end{align}

Now, we estimate  these terms on the right-hand side of \eqref{4.16} one by one.  First, it follows from \eqref{3.3}, \eqref{4.2}, \eqref{4.3}, \eqref{4.4} and \eqref{4.6} that
\begin{align}\la{4.17}
&\left|\int\left(2\n\uu\cdot\na\uu_t\cdot\uu_{tt}+\n\uu_t\cdot\na \uu \cdot\uu_{tt}\right)d\xx\right|+\left|\int\n\uu\cdot\na(\uu\cdot\na\uu) \cdot \uu_{tt}d\xx\right|\\
\leq &\ve\|\sqrt\n \uu_{tt}\|_{L^2}^2 +C(\ve)\left(\|\sqrt\n\uu\| _{L^\infty}^2 \|\na\uu_t\|_{L^2}^2+\|\sqrt\n\uu_t\|^2_{L^4}\|\na\uu\| _{L^4}^2\right)\no
&+C(\ve)\left(\|\sqrt\n\uu\|_{L^\infty}^2\|\na\uu\|_{L^4}^2 +\|\n^{1/4}\uu\|_{L^\infty}^4\|\na^2\uu\|_{L^2}^2\right) \no
\leq &\ve\|\sqrt\n \uu_{tt}\|_{L^2}^2 +C(\ve)\left(1+\|\na\uu_t\|_{L^2}^2 \right).\nonumber
\end{align}

Then one yields that
\begin{align}\la{4.18}
&-\int\n\uu\cdot\na\uu_{tt}\cdot\uu_td\xx-\int\n\uu\cdot\na\uu_{tt} \cdot(\uu \cdot\na)\uu d\xx\\
=&-\frac{d}{dt}\int\left( \n\uu\cdot\na\uu_{t}\cdot\uu_t+\n\uu \cdot\na\uu_{t}\cdot(\uu \cdot\na)\uu\right)d\xx +\int(\n\uu)_t\cdot\na\uu_{t}\cdot\uu_td\xx\no
&+\int(\n\uu)_t\cdot\na\uu_{t}\cdot(\uu \cdot\na)\uu d\xx+\int\n\uu\cdot\na\uu_{t}\cdot\uu_{tt}d\xx\no
&+\int\n\uu\cdot\na\uu_{t}\cdot(\uu_t \cdot\na)\uu d\xx +\int\n\uu\cdot\na\uu_{t}\cdot(\uu \cdot\na)\uu_t d\xx.\nonumber
\end{align}
Noting that those terms without time derivatives on the right hand side in \eqref{4.18}, H\"{o}lder inequality together with \eqref{3.3} and \eqref{4.4}--\eqref{4.6} implies
\begin{align}\la{4.19}
\int(\n\uu)_t\cdot\na\uu_{t}\cdot\uu_td\xx
\leq & C\|\n\bar\xx^a\|_{L^\infty}\|\bar\xx^{-a/2}\uu_t\|_{L^4}^2 \|\na\uu_t\|_{L^2}\\
&+C\|\bar\xx^{(1+a)/2}\n_t\|_{L^2} \|\bar\xx^{-1/2}\uu\|_{L^\infty}\|\bar\xx^{-a/2}\uu_t\|_{L^4} \|\na\uu_t\|_{L^4}\no
\leq &\delta\|\na^2\uu_t\|_{L^2}^2+C(\delta)\|\na\uu_t\|_{L^2}^4 +C(\delta),\nonumber
\end{align}
and
\begin{align}\la{4.20}
&\int(\n\uu)_t\cdot\na\uu_{t}\cdot(\uu \cdot\na)\uu d\xx\\
\leq  & C\|\n\bar\xx^a\|_{L^\infty} \|\bar\xx^{-a/2}\uu_t\|_{L^4} \|\bar\xx^{-a/2}\uu\cdot\na\uu\|_{L^4}\|\na\uu_t\|_{L^2}\no
&+C\|\bar\xx^{(1+a)/2}\n_t\|_{L^2}\|\bar\xx^{-1/2}\uu\|_{L^\infty} \|\bar\xx^{-a/2}\uu\|_{L^\infty}\|\na\uu\|_{L^4}\|\na\uu_t\|_{L^4}\no
\leq &\delta\|\na^2\uu_t\|_{L^2}^2+C(\delta)\|\na\uu_t\|_{L^2}^4 +C(\delta).\nonumber
\end{align}
It is easy from\eqref{4.4} to lead to
\begin{align}\la{4.21}
\int\n\uu\cdot\na\uu_{t}\cdot\uu_{tt}d\xx\leq & C\|\sqrt\n \uu_{tt}\|_{L^2}\|\sqrt\n\uu\|_{L^\infty}\|\na\uu_t\|_{L^2}\\
\leq &\ve\|\sqrt\n \uu_{tt}\|_{L^2}^2 +C(\ve)\|\na\uu_t\|_{L^2}^2.\nonumber
\end{align}
Next, it follows from \eqref{4.3}, \eqref{4.4} and \eqref{4.6} that
\begin{align}\la{4.22}
\int\n\uu\cdot\na\uu_{t}\cdot(\uu_t \cdot\na)\uu d\xx\leq &C\|\n\bar\xx^a\|_{L^\infty}\|\bar\xx^{-a/2} \uu\|_{L^\infty} \|\bar\xx^{-a/2} \uu_t\|_{L^4}\|\na\uu\|_{L^4}\|\na\uu_t\|_{L^2}\no
\leq & C+C(\ve)\|\na\uu_t\|_{L^2}^2,
\end{align}
and similarly,
\begin{align}\la{4.23}
\int\n\uu\cdot\na\uu_{t}\cdot(\uu \cdot\na)\uu_t d\xx \leq &C\|\n\bar\xx^a\|_{L^\infty}\|\bar\xx^{-a/2} \uu\|_{L^\infty}^2\|\na\uu_t\|_{L^2}^2\\
\leq & C\|\na\uu_t\|_{L^2}^2.\nonumber
\end{align}
Inserting \eqref{4.19}--\eqref{4.23} into \eqref{4.18}, we conclude that
\begin{align}\la{4.24}
&-\int\n\uu\cdot\na\uu_{tt}\cdot\uu_td\xx-\int\n\uu\cdot\na\uu_{tt} \cdot(\uu \cdot\na)\uu d\xx\\
=&-\frac{d}{dt}\int\left( \n\uu\cdot\na\uu_{t}\cdot\uu_t+\n\uu \cdot\na\uu_{t}\cdot(\uu \cdot\na)\uu\right)d\xx+\ve\|\sqrt\n \uu_{tt}\|_{L^2}^2\no
&+C(\ve, \delta)\|\na\uu_t\|_{L^2}^4+\delta\|\na^2\uu_t\|_{L^2}^2 +C(\ve, \delta).\nonumber
\end{align}

On the other hand, from \eqref{2.2}$_3$, \eqref{4.1}, \eqref{4.3} and \eqref{4.5} , we obtain that
\begin{align}\la{4.26}
&-\int\hh_t\cdot\na\uu_{tt}\cdot \hh d\xx-\int\hh\cdot\na\uu_{tt}\cdot\hh_td\xx\\
=&-\frac{d}{dt}\int\left(\hh_t\cdot\na\uu _t\cdot\hh +\hh\cdot\na\uu_t\cdot\hh_t\right)d\xx\no
&+\int\hh_{tt}\cdot\na\uu_t\cdot\hh d\xx+2\int\hh_t\cdot\na\uu_t\cdot\hh_t d\xx+\int\hh\cdot\na\uu_t\cdot \hh_{tt}d\xx\no
\leq &-\frac{d}{dt}\int\left(\hh_t\cdot\na\uu _t\cdot\hh +\hh\cdot\na\uu_t\cdot\hh_t\right)d\xx\no
&+C\|\hh_t\bar\xx^{(a+1)/2}\|_{L^q} \|\na\uu_t\|_{L^2}\left(\|\uu\bar\xx^{-1/2}\|_{L^\infty} \|\na\hh\|_{(q-2)/2q}+\|\hh\|_{L^\infty}\|\na\uu\|_{(q-2)/2q}\right)\no
&+C\|\bar\xx^{a/2}\hh\|_{L^\infty}\|\bar\xx^{-a/2}\uu_t\|_{L^4} \|\na\hh\|_{L^4} \|\na\uu_t\|_{L^2}+C\|\hh\|_{L^\infty}^2\|\na \uu_t\|_{L^2}^2\no
&+C\|\bar\xx^{a/2}\hh\|_{L^\infty} \|\bar\xx^{-a/2}\uu\|_{L^\infty} \|\na\hh_t\|_{L^2} \|\na\uu_t\|_{L^2}\no
&+C\|\hh\|_{L^\infty}\|\hh_t \bar\xx^{a/2}\|_{L^q}\|\na\uu_t\|_{L^2}\|\na\uu\|_{(q-2)/2q}\no
\leq &-\frac{d}{dt}\int\left(\hh_t\cdot\na\uu _t\cdot\hh +\hh\cdot\na\uu_t\cdot\hh_t\right)d\xx+C \left(\|\na\uu_t\|_{L^2}^2+1\right),\nonumber
\end{align}
where in the last inequality we have used \eqref{4.6} and the following simple fact,
\begin{align}\la{4.27}
\|\na\hh_t\|_{L^2}\leq& C\||\na\uu||\na\hh|\|_{L^2}+C\||\uu||\na^2\hh|\| _{L^2}+C\||\hh||\na^2\uu|\|_{L^2}\\
\leq &C\|\na\uu\|_{L^{2q/(q-2)}}\|\na\hh\|_{L^q}+C\|\bar\xx^{-\delta_0} \uu\|_{L^\infty} \|\bar\xx^{\delta_0}\na^2\hh\|_{L^2}+C\|\hh\|_{L^\infty} \|\na^2\uu\|_{L^2}\no
\leq &C,\nonumber
\end{align}
thanks to \eqref{4.1}, \eqref{4.3}, \eqref{4.4} and \eqref{4.7}.

Substituting \eqref{4.17} and \eqref{4.24}--\eqref{4.26} into \eqref{4.16} and choosing $\ve$ suitably small, we can yield that
\be\la{4.28}
\Phi'(t)+\|\sqrt\n \uu_{tt}\|_{L^2}^2\leq C\delta\|\na^2\uu_t\|_{L^2}^2 +C\|\na\uu_t\|_{L^2}^4+C,
\ee
where
\begin{align*}
\Phi(t)\triangleq&\mu\|\na\uu_t\|_{L^2}^2+\int\left( \n\uu\cdot\na\uu_{t}\cdot\uu_t+\n\uu \cdot\na\uu_{t}\cdot(\uu \cdot\na)\uu\right)d\xx\no
&+\int\left(\hh_t\cdot\na\uu _t\cdot\hh +\hh\cdot\na\uu_t\cdot\hh_t\right)d\xx
\end{align*}
satisfies
\be\la{4.29}
C(\mu)\|\na\uu_t\|_{L^2}^2-C\leq \Phi(t)\leq C\|\na\uu_t\|_{L^2}^2 +C,
\ee
from the following estimate which has been yielded  from \eqref{4.2}--\eqref{4.5}
\begin{align*}
&\left|\int\left( \n\uu\cdot\na\uu_{t}\cdot\uu_t+\n\uu \cdot\na\uu_{t}\cdot(\uu \cdot\na)\uu\right)d\xx\right|+\left|\int\left(\hh_t\cdot\na\uu _t\cdot\hh +\hh\cdot\na\uu_t\cdot\hh_t\right)d\xx\right|\no
\leq &C\|\sqrt\n \uu\|_{L^\infty}\|\na\uu_t\|_{L^2}\|\sqrt\n\uu_t\|_{L^2}
+C\|\sqrt\n \uu\|_{L^\infty}^2\|\na\uu_t\|_{L^2}\|\na\uu\|_{L^2}\no
&+C\|\hh\|_{L^\infty}\|\hh_t\|_{L^2}\|\na\uu_t\|_{L^2}
\leq \ve\|\na\uu_t\|_{L^2}^2+C(\ve).
\end{align*}

To end this proof, it remains to estimate the first term on the right-hand side of \eqref{4.28}. In fact,  rewrite the equation \eqref{3.25} as follows,
\begin{align*}
-\mu\triangle\uu_t+\nabla \pi_t&=-\n\uu_{tt}-\n\uu\cdot\na\uu_t \\
&-\n_t(\uu_t+\uu\cdot\na\uu)-\n\uu_t\cdot\na\uu+\hh_t\cdot\na\hh +\hh\cdot\na\hh_t.\nonumber
\end{align*}

By Lemma \ref{l2.5}  we obtain that
\begin{align}\la{4.30}
&\|\na^2\uu_t\|_{L^2}^2\\ \leq&C(\|-\n\uu_{tt}-\n\uu\cdot\na\uu_t-\n_t(\uu_t+\uu\cdot\na\uu) -\n\uu_t\cdot\na\uu+\hh_t\cdot\na\hh+\hh\cdot\na\hh_t\|_{L^2}^2)\no
\leq & C\|\n\|_{L^\infty} \|\sqrt\n\uu_{tt} \|_{L^2}^2 +C\|\n\| _{L^\infty} \|\sqrt\n\uu \|_{L^\infty}^2 \|\na\uu_t\|_{L^2}^2 +C\|\bar\xx^{(a+1)/2}\n_t\|_{L^q}^2\|\bar\xx^{-1} \uu_t\|_{L^{2q/(q-2)}}^2 \no &+C\|\bar\xx^{(a+1)/2}\n_t\|_{L^q}^2\|\bar\xx^{-1}\uu\|_{L^\infty}^2  \|\na\uu\|_{L^{2q/(q-2)}}^2+C\|\bar\xx^{(a+1)/2}\n\|_{L^q}^2 \|\bar\xx^{-1}\uu\|_{L^\infty}^2\|\na\uu\|_{L^{2q/(q-2)}}^2\no
&+C\|\sqrt\n\|_{L^\infty}\|\sqrt\n\uu_t\|_{L^2}^2\|\nabla\uu\|_{L^2}^2 +C\|\bar\xx^{(a+1)/2}\hh_t\|_{L^q}^2 \|\na\hh\|_{L^{2q/(q-2)}}^2+C\|\hh\|_{L^\infty}^2\|\na\hh_t\|_{L^2}^2\no
\leq & C\|\sqrt\n\uu_{tt} \|_{L^2}^2+C\|\na\uu_t\|_{L^2}^4+C,\nonumber
\end{align}
where we have used Gagliardo-Nirenberg inequality, \eqref{4.1}-- \eqref{4.7} and \eqref{4.27}.

Inserting \eqref{4.30} into \eqref{4.28} and choosing $\delta$ suitably small  lead to
\be\la{4.31}
\Phi'(t)+\|\sqrt\n \uu_{tt}\|_{L^2}^2\leq C\|\na\uu_t\|_{L^2}^4+C.
\ee

Multiplying \eqref{4.31} by $t$ and integrating the resultant inequality over $(0, T_0)$, we obtain from Gronwall inequality, \eqref{4.2} and \eqref{4.29}
$$
\sup_{0\leq t\leq T_0}t\|\na\uu_t\|_{L^2}^2+\int_0^{T_0}t\|\sqrt\n \uu_{tt}\|_{L^2}^2dt\leq C.
$$
Combining \eqref{4.30},  we  yield \eqref{4.15} and the complete the proof of Lemma \ref{l4.2}. \hfill $\Box$

\bl\la{l4.3}
It holds that
\be\la{4.32}
\sup_{0\leq t\leq T_0}\left( \|\na^2\n\|_{L^q}+\|\na^2\hh\|_{L^q}\right)\leq C.
\ee
\el
\pf Applying the differential operator $\na^2$ to \eqref{a4.7} and \eqref{4.8}, respectively, and multiplying each equality by $q|\na^2\bar\n| ^{q-2} \na^2\bar\n$ and $q|\na^2\bar\hh|^{q-2}\na^2\bar\hh$, and integrating the resultant equalities over $B_R$ lead to
\begin{align}\la{4.33}
&\frac{d}{dt}\left(\|\na^2\bar\hh\|_{L^q}+\|\na^2\bar\n\|_{L^q}\right)\\
\leq &C\|\na\uu\|_{L^\infty}\left(\|\na^2\bar\hh\|_{L^q} +\|\na^2 \bar\n\|_{L^q}\right)+ \left( \|\na\bar\hh\| _{L^\infty} +\|\na \bar\n\|_{L^\infty}\right)\|\na^2\uu\|_{L^q}\no
\leq &C\left(1+\|\na^2\uu\|_{L^q}\right)\left(1+\|\na^2\bar\hh\|_{L^q} +\|\na^2 \bar\n\|_{L^q}\right)+C\|\na^3\uu \|_{L^q}.\nonumber
\end{align}
Due to \eqref{2.7}, the last term on the right-hand side of \eqref{4.33} can be estimated as follows:
\begin{align}\la{4.34}
\|\na^3\uu\|_{L^q}+\|\nabla^2\pi\|_{L^q}\leq &C\|\na(\n\uu_t)\|_{L^q}+C\|\na(\n\uu\cdot\na \uu)\|_{L^q}+C\|\na (\hh\cdot\na\hh )\|_{L^q}\\
\leq & C\|\bar\xx^{-a}\uu_t\|_{L^\infty}\|\bar\xx^a\na\n\|_{L^q} +C\|\bar\xx^a\n\|_{L^\infty}\|\bar\xx^{-a}\na\uu_t\|_{L^q}\no
&+C\|\bar\xx^a\na\n\|_{L^q}\|\bar\xx^{-a}\uu\|_{L^\infty} \|\na\uu\|_{L^\infty}+C\|\n\|_{L^\infty}\|\na\uu\|_{L^{2q}}^2 \no &+C\|\bar\xx^a\n\|_{L^\infty}\|\bar\xx^{-a}\uu\|_{L^\infty} \|\na^2\uu\|_{L^q}\no
&+C\|\bar\xx^{a}\hh\|_{L^\infty}\|\bar\xx^{-a}\na^2\hh\|_{L^q} +C\|\na\hh\|_{L^{2q}}^2\no
\leq &C\|\na\uu_t\|_{L^q}+C\|\bar\xx^{-a}\uu_t\|_{L^q} +\frac12\|\na^3\uu\|_{L^q}+C\|\na^2\bar\hh\|_{L^q}+C\no
\leq &C\|\na\uu_t\|_{L^2}^{2/q}\|\na^2\uu_t\|_{L^2}^{(q-2)/q} +C\|\na\uu_t\|_{L^2}+\frac12\|\na^3\uu\|_{L^q}\no
&+C\|\na^2\bar\hh\|_{L^q}+C,\nonumber
\end{align}
where  \eqref{3.3}, \eqref{4.3}, \eqref{4.4}, \eqref{4.6} and \eqref{4.7} were used.

While, it follows from \eqref{4.15} that
\begin{align}\la{4.35}
&\int_0^{T_0}\left(\|\na\uu_t\|_{L^2}^{2/q}\|\na^2\uu_t\|_{L^2} ^{(q-2)/q}\right)^{(q+1)/q}dt\\
\leq &C\sup_{0\leq t\leq T_0}\left(t\|\na\uu_t\|_{L^2}^2\right) ^{(q+1)/q^2}\int_0^{T_0}\left(t\|\na^2\uu_t\|_{L^2}^2 +t^{-(q^2+q)/(q^2+q+2)}\right)dt\no
\leq &C.\nonumber
\end{align}
Putting \eqref{4.34} into \eqref{4.33}, we get \eqref{4.32} from Gronwall inequality, \eqref{3.3}, \eqref{4.2} and \eqref{4.35}. Therefore, the proof of Lemma \ref{l4.3} is finished. \hfill $\Box$

\bl\la{l4.4}
It holds that
\begin{align}\la{4.36}
&\sup_{0\leq t\leq T_0}t\left( \|\na^3\uu\|_{L^2\cap L^q}+\|\na\uu_t\|_{H^1} +\|\na^2(\n\uu)\|_{L^{(q+2)/2}}\right)\\
&\qquad\qquad\quad+\int_0^{T_0}t^2\left(\|\na\uu_{tt}\| _{L^2}^2+\|\bar\xx^{-1}\uu_{tt}\|_{L^2}^2 \right)dt \leq C.\nonumber
\end{align}
\el
\pf First, we claim that
\be\la{4.37}
\sup_{0\leq t\leq T_0}t^2\|\sqrt\n\uu_{tt}\|_{L^2}^2+\int_0^{T_0} t^2\|\na\uu_{tt}\|_{L^2}^2 dt\leq C,
\ee
which together with \eqref{2.5}, \eqref{4.15} and \eqref{4.30} yields that
\be\la{4.38}
\sup_{0\leq t\leq T_0}t\|\na\uu_t\|_{H^1}+\int_0^{T_0}t^2\|\bar\xx^{-1} \uu_{tt}\|_{L^2}^2dt \leq C.
\ee
This combined with \eqref{4.13}, \eqref{4.32}, \eqref{4.34} and \eqref{4.35} leads to
\be\la{4.39}
\sup_{0\leq t\leq T_0}t\|\na^3\uu\|_{L^2\cap L^q}\leq C,
\ee
which together with \eqref{3.3}, \eqref{4.1} and \eqref{4.32}, shows
\begin{align}\la{4.40}
t\|\na^2(\n\uu)\|_{L^{(q+2)/2}}\leq &C t\||\na^2\n||\uu|\|_{L^{(q+2)/2}} +Ct\||\na\n||\na\uu|\|_{L^{(q+2)/2}}+Ct\|\n|\na^2\uu|\|_{L^{(q+2)/2}}\no
\leq &Ct\|\bar\xx^{\delta_0}\na^2\n\|_{L^2}^{2/(q+2)}\|\na^2\n\|_{L^q} ^{q/(q+2)}\|\bar\xx^{-2\delta_0/(q+2)}\uu\|_{L^\infty}\\
&+Ct\|\na\n\|_{L^q}\|\na\uu\|_{L^{q(q+2)/(q-2)}}+Ct\|\na^2\uu\| _{L^{(q+2)/2}}\no
\leq &C.\nonumber
\end{align}
Therefore, we complete the proof of \eqref{4.36} from \eqref{4.37}--\eqref{4.40}.

Now, we focus on the estimates of \eqref{4.37}. In fact, differentiating \eqref{3.25} with respect to $t$ yields that
\begin{align*}
&\n\uu_{ttt}+\n\uu\cdot\na\uu_{tt}-\mu\triangle\uu_{tt}+\nabla\pi_{tt}\no
=&2\div(\n\uu)\uu_{tt}+\div(\n\uu)_t\uu_t-2(\n\uu)_t\cdot\na\uu_t -\n_{tt}\uu\cdot\na\uu-2\n_t\uu_t\cdot\na\uu\no
&-\n\uu_{tt}\cdot\na\uu+\hh_{tt}\cdot\na\hh+2\hh_t\cdot\na\hh_t+\hh \cdot\na\hh_{tt},
\end{align*}
which multiplied by $\uu_{tt}$ and integrated by parts over $B_R$, shows that
\begin{align}\la{4.41}
&\frac12\frac{d}{dt}\int\n|\uu_{tt}|^2d\xx+\int\left(\mu|\na\uu_{tt}|^2 \right)d\xx\\
=&-4\int\n\uu\cdot\na\uu_{tt}\cdot\uu_{tt}d\xx-\int(\n\uu)_t\cdot\left( \na(\uu_t\cdot\uu_{tt})+2\na\uu_t\cdot\uu_{tt}\right)\no
&-\int(\n\uu)_t\cdot\na(\uu\cdot\na\uu\cdot\uu_{tt})d\xx-2\int\n_t\uu_t \cdot \na\uu\cdot\uu_{tt}d\xx-\int\n\uu_{tt}\cdot\na\uu\cdot\uu_{tt}d\xx\no
&+\int\hh_{tt}\cdot\na\hh\cdot\uu_{tt}d\xx+2\int\hh_t\cdot\na\hh_t\cdot \uu_{tt}d\xx+\int\hh\cdot\na\hh_{tt}\cdot\uu_{tt}d\xx\no
\triangleq&\sum_{i=1}^{8}I_i.\nonumber
\end{align}

Now, we  estimate each term on the right-hand side of \eqref{4.41} as follows:

First, it follows from \eqref{4.4} that
\be\la{4.42}
|I_1|\leq C\|\sqrt\n\uu\|_{L^\infty}\|\sqrt\n\uu_{tt}\|_{L^2}\|\na \uu_{tt}\|_{L^2}\leq \ve\|\na \uu_{tt}\|_{L^2}^2+C(\ve) \|\sqrt\n\uu_{tt}\|_{L^2}^2.
\ee

It is easy to check that from \eqref{2.6}, \eqref{3.9} and \eqref{4.5}
\begin{align}\la{4.43}
|I_2|\leq &C\|\bar\xx(\n\uu)_t\|_{L^q}\left(\|\bar\xx^{-1}\uu_{tt}\| _{2q/(q-2)}\|\na\uu_t\|_{L^2}+\|\bar\xx^{-1}\uu_{t}\| _{2q/(q-2)}\|\na\uu_{tt}\|_{L^2}\right)\\
\leq &C\left(1+\|\na\uu_t\|_{L^2}^2\right)\left(\|\sqrt\n\uu_{tt}\| _{L^2}+\|\na\uu_{tt}\|_{L^2} \right)\no
\leq &\ve\left(\|\sqrt\n\uu_{tt}\| _{L^2}^2+\|\na\uu_{tt}\|_{L^2}^2 \right)+C(\ve)\left(1+\|\na\uu_t\|_{L^2}^4\right),\nonumber
\end{align}
where we have  also used  the following facts:
\begin{align}\la{4.44}
\|\bar\xx(\n\uu)_t\|_{L^q}\leq &C\|\bar\xx|\n_t||\uu|\|_{L^q} +C\|\bar\xx \n|\uu_t|\|_{L^q}\\
\leq &C\|\n_t\bar\xx^{(1+a)/2}\|_{L^q}\|\bar\xx^{-(a-1)/2}\uu\| _{L^\infty} +C\|\n\bar\xx^a\|_{L^{2q/(3-\ti a)}}\|\uu_t\bar\xx^{1-a}\|_{L^{2q/(\ti a-1)}}\no
\leq &C+C\|\na\uu_t\|_{L^2}.\nonumber
\end{align}
by \eqref{4.4}--\eqref{4.6}, where $\ti a=\min\{2, a\}$.

Then, it follows from \eqref{3.9}, \eqref{4.3}, \eqref{4.4} and \eqref{4.44} that
\begin{align}\la{4.45}
|I_3|\leq &C\int|(\n\uu)_t|\left(|\uu||\na^2\uu||\uu_{tt}|+|\uu||\na\uu| |\na\uu_{tt}|+|\na\uu|^2|\uu_{tt}|\right)d\xx\\
\leq &C\|\bar\xx(\n\uu)_t\|_{L^q}\|\bar\xx^{-1/q}\uu\|_{L^\infty}\|\na^2\uu\| _{L^2} \|\bar\xx^{-(q-1)/q}\uu_{tt}\|_{L^{2q/(q-2)}}\no
&+C\|\bar\xx(\n\uu)_t\|_{L^q}\|\bar\xx^{-1}\uu\|_{L^\infty}\|\na\uu\| _{L^{2q/ (q-2)}}\|\na\uu_{tt}\|_{L^2}\no
&+C\|\bar\xx(\n\uu)_t\|_{L^q}\|\na\uu\|_{L^4}^2\|\bar\xx^{-1}\uu_{tt}\| _{L^{2q/(q-2)}}\no
\leq &C(1+\|\na\uu_t\|_{L^2})\left(\|\sqrt\n\uu_{tt}\| _{L^2}+\|\na\uu_{tt}\|_{L^2} \right)\no
\leq &\ve\left(\|\sqrt\n\uu_{tt}\| _{L^2}^2+\|\na\uu_{tt}\|_{L^2}^2 \right)+C(\ve)\left(1+\|\na\uu_t\|_{L^2}^2\right).\nonumber
\end{align}

Clearly, it follows from Cauchy inequality, together with \eqref{3.9}, \eqref{4.5} and \eqref{4.6} that
\begin{align}\la{4.46}
|I_4|\leq &C\int|\n_t||\uu_t||\na\uu||\uu_{tt}|d\xx\\
\leq & C\|\bar\xx\n_t\|_{L^q}\|\bar\xx^{-1/2}\uu_t\|_{L^{4q/(q-2)}} \|\na\uu\|_{L^2}\|\bar\xx^{-1/2}\uu_{tt}\|_{L^{4q/(q-2)}}\no
\leq & C\left(1+\|\na\uu_t\|_{L^2}\right)\left(\|\sqrt\n\uu_{tt}\| _{L^2}+\|\na\uu_{tt}\|_{L^2} \right)\no
\leq &\ve\left(\|\sqrt\n\uu_{tt}\| _{L^2}^2+\|\na\uu_{tt}\|_{L^2}^2 \right)+C(\ve)\left(1+\|\na\uu_t\|_{L^2}^2\right).\nonumber
\end{align}

Then, Gagliardo-Nirenberg inequality together with \eqref{4.3} gives
\begin{align}\la{4.47}
|I_5|\leq C\|\na\uu\|_{L^\infty}\|\sqrt\n\uu_{tt}\|_{L^2}^2\leq C\left(1+\|\na^2\uu\|_{L^q}\right)\|\sqrt\n\uu_{tt}\|_{L^2}^2.
\end{align}

Finally, it follows from \eqref{2.2}$_3$, \eqref{3.3} and \eqref{4.3}--\eqref{4.6} that
\begin{align}\la{4.49}
\|\hh_{tt}\|_{L^2}\leq & C\|\bar\xx^{-\theta_0}\uu_t\| _{2q/[(q-2)\theta_0]} \|\bar\xx^{\theta_0} \na \hh\|_{L^{2q/[q-(q-2)\theta_0]}}+C\|\na \uu_t\|_{L^2}\\
&+C\|\bar\xx^{-\delta_0/2}\uu\|_{L^\infty}\|\bar\xx^{\delta_0/2}\na \hh_t\|_{L^2}+C\|\hh_t\|_{L^q}\|\na\uu\|_{L^{2q/(q-2)}}\no
\leq &C\left(1+\|\na\uu_t\|_{L^2}\right),\nonumber
\end{align}
where in the last inequality we have used the following fact:
\begin{align*}
\|\bar\xx^{\delta_0/2}\na\hh_t\|_{L^2}\leq& C\|\bar\xx^{\delta_0/2}|\na\uu||\na\hh|\|_{L^2}+C\|\bar\xx^{\delta_0/2} |\uu||\na^2\hh|\| _{L^2}+C\|\bar\xx^{\delta_0/2} |\hh||\na^2\uu|\|_{L^2}\no
\leq &C\|\na\uu\|_{L^{2q/(q-2)}}\|\bar\xx^{a}\na\hh\|_{L^q} +C\|\bar\xx^{-\delta_0/2} \uu\|_{L^\infty} \|\bar\xx^{\delta_0} \na^2\hh\|_{L^2}\no
&+C\|\bar\xx^{\delta_0/2}\hh\|_{L^\infty} \|\na^2\uu\|_{L^2}\no
\leq &C,
\end{align*}
due to \eqref{4.1}, \eqref{4.3}, \eqref{4.4} and \eqref{4.7}. Then \eqref{4.49} and integration by parts lead to
\begin{align}\la{4.50}
|I_6|+|I_7|+|I_8|\leq &C\|\hh_t\|_{L^4}^2\|\na\uu_{tt}\|_{L^2} +C\|\hh\|_{L^\infty}\|\hh_{tt}\|_{L^2}\|\na\uu_{tt}\|_{L^2}\\
\leq &C\left(\|\na\uu_t\|_{L^2}+\|\hh_t\|_{L^2}^2 +\|\na\hh_t\|_{L^2}^2\right) \|\na\uu_{tt}\|_{L^2}\no
\leq &\ve\|\na\uu_{tt}\|_{L^2}^2 +C(\ve)\left(1+\|\na\uu_t\|_{L^2}^2\right),\nonumber
\end{align}
in terms of \eqref{4.5}, \eqref{4.7}, \eqref{4.27} and \eqref{4.49}.

Substituting \eqref{4.42}, \eqref{4.43}, \eqref{4.45}--\eqref{4.47} and \eqref{4.50} into \eqref{4.41}, choosing $\ve$ small enough, and multiplying the resultant inequality by $t^2$, we get \eqref{4.37} after using Gronwall inequality and \eqref{4.15}. Therefore, the proof of Lemma \ref{l4.4} is completed. \hfill $\Box$

\section{Proofs of main theorems }

With all the a priori estimates obtained in Section 3 and 4 at hand, now we are ready to prove the main results of this paper in this section.

{\it Proof of Theorem \ref{thm1}.}
Let $(\n_0, \uu_0, \hh_0)$ be as in Theorem \ref{thm1}. Without loss of generality, we assume that the initial density $\n_0$ satisfies
$$
\int_{\mr^2}\n_0 d\xx=1,
$$
which means that there has a positive constant $N_0$ such that
\be\la{5.1}
\int_{B_{N_0}}\n_0d\xx\geq\frac34\int_{\mr^2}\n_0 d\xx=\frac34.
\ee
We construct that $\n_0^R=\hat\n_0^R+R^{-1}e^{-|\xx|^2}$ where $0\leq \hat\n_0^R\in C_0^\infty(\mr^2)$ satisfies
\be\la{5.2}
\left\{
\begin{array}{lll}
\int_{B_{N_0}}\hat\n_0^R d\xx\geq 1/2,\\
\bar\xx^a\hat\n_0^R\rightarrow\bar\xx^a\n_0 \quad {\rm in\ } L^1(\mr^2)\cap H^1(\mr^2)\cap W^{1, q}(\mr^2),
\end{array}
\right.\quad{\rm as}\ R\rightarrow\infty.
\ee
Then, we choose $\hh_0^R\in \{\mathbf{w}\in C_0^\infty(B_R)|\div \mathbf{w}=0\}$ satisfying
\be\la{5.3}
\hh_0^R\bar\xx^a\rightarrow \hh_0\bar\xx^a \quad  {\rm in \ }H^1(\mr^2)\cap W^{1, q}(\mr^2),\qquad {\rm as}\ R \rightarrow\infty.
\ee
Next, since $\na\uu_0\in L^2(\mr^2)$, choosing $\mathbf v_i^R\in C_0^\infty(B_R)$ $(i=1, 2)$ such that
\be\la{5.4}
\lim_{R\rightarrow\infty}\|\mathbf v_i^R-\p_i\uu_0\|_{L^2(\mr^2)}=0,\quad i=1, 2,
\ee
we consider the unique smooth solution $\uu_0^R$ of the following Stokes problem:
\be\la{5.5}
\left\{
\begin{array}{lll}
-\triangle \uu_0^R+\nn\uu_0^R+\nabla\pi^R_0=\sqrt{\nn}\mathbf{h}^R-\p_i\mathbf v_i^R, &\ {\rm in\ }B_R,\\
\div\uu=0,&\ {\rm in\ }B_R,\\
\uu_0^R=0, &\ {\rm on\ }\p B_R,
\end{array}
\right.
\ee
where $\mathbf{h}^R=(\sqrt{\n_0}\uu_0)*j_{1/R}$ with $j_\delta$ being the standard mollifying kernel of width $\delta$. Extending $\uu_0^R$ to $\mr^2$ by defining $0$ outside of $B_R$ and denoting it by $\ti\uu_0^R$, in the \cite{lv}, we have known that
\be\la{5.6}
\lim_{R\rightarrow\infty}\left(\|\na(\ti\uu_0^R-\uu_0)\|_{L^2(\mr^2)} +\|\sqrt{\nn}\ti\uu_0^R-\sqrt{\n_0}\uu_0\|_{L^2(\mr^2)} \right)=0.
\ee

In the sequence, thanks to Lemma \ref{l2.1}, the initial-boundary value problem \eqref{2.2} with the initial data $(\n_0^R, \uu_0^R, \hh_0^R)$ has a classical solution $(\n^R, \uu^R, \hh^R)$ on $B_R\times[0, T_R]$. Moreover, Proposition \ref{prop} gives that there has a $T_0$ independent of $R$ such that \eqref{3.3} holds for $(\n^R, \uu^R, \hh^R)$. Extending $(\n^R, \uu^R, \hh^R)$ by zero on $\mr^2\backslash B_R$ and denoting it by
$$
\ti\n^R\triangleq\varphi_R\n^R, \quad \ti\uu^R, \quad \ti\hh^R\triangleq\varphi_R\hh^R,
$$
with $\varphi_R$ as in \eqref{3.6}, we first deduce from \eqref{3.3} that
\begin{align}\la{5.11}
&\sup_{0\leq t\leq T_0}\left(\|\sqrt{\ti\n^R}\ti\uu^R\|_{L^2} +\|\na \ti\uu^R\|_{L^2}\right)\\
\leq &C+C\sup_{0\leq t\leq T_0}\|\na\uu^R\|_{L^2}\no
\leq &C,\nonumber
\end{align}
and
\be\la{5.12}
\sup_{0\leq t\leq T_0}\|\ti\n^R\bar\xx^a\|_{L^1\cap L^\infty}\leq C.
\ee

Next, for $p\in [2, q]$, it follows from \eqref{3.3} and \eqref{3.33} that
\begin{align}\la{5.13}
&\sup_{0\leq t\leq T_0}\left(\|\na(\bar\xx^a\ti\n^R)\|_{L^p(\mr^2)}+ \|\bar\xx^a\na\ti\n^R\|_{L^p(\mr^2)}+\|\na(\bar\xx^a\ti\hh^R)\| _{L^p(\mr^2)} +\|\bar\xx^a\na\ti\hh^R\|_{L^p(\mr^2)}\right)\no
\leq &C\sup_{0\leq t\leq T_0}\left(\|\bar\xx^a\na\n^R\|_{L^p(\O)} +\|\bar\xx^a\n^R\na\varphi_R\|_{L^p(\O)}+\|\n^R\na\bar\xx^a\|_{L^p(\O)} \right)\\
&+C\sup_{0\leq t\leq T_0}\left(\|\bar\xx^a\na\hh^R\|_{L^p(\O)} +\|\bar\xx^a\hh^R\na\varphi_R\|_{L^p(\O)}+\|\hh^R\na\bar\xx^a\|_{L^p(\O)} \right)\no
\leq &C +C\|\bar\xx^a\n^R\|_{L^p(\O)}+C\|\bar\xx^a\hh^R\|_{L^p(\O)}\no
\leq &C.\nonumber
\end{align}

Then, it follows from \eqref{3.3} and \eqref{3.35} that
\be\la{5.14}
\int_0^{T_0}\left(\|\na^2 \ti\uu^R\|_{L^q(\mr^2)}^{(q+1)/q}+t \|\na^2\ti\uu^R\|_{L^q(\mr^2)}^2+\|\na^2\ti\uu^R\|_{L^2(\mr^2)} ^2\right)dt\leq C,
\ee
and that for $p\in [2, q]$,
\begin{align}\la{5.15}
&\int_0^{T_0}\left(\|\bar\xx\ti\n_t^R\|_{L^p(\mr^2)}^2 +\|\bar\xx\ti\hh_t^R\|_{L^p(\mr^2)}^2\right)dt\\
\leq & C\int_0^{T_0}\left(\||\bar\xx||\uu^R||\na\n^R|\|_{L^p(B_R)}^2 +\|\bar\xx\n^R\div\uu^R\|_{L^p(B_R)}^2\right)dt\no
&+C\int_0^{T_0}\left(\||\bar\xx||\uu^R||\na\hh^R|\|_{L^p(B_R)}^2 +\|\bar\xx|\hh^R||\na\uu^R|\|_{L^p(B_R)}^2\right)dt\no
\leq &C\int_0^{T_0}\|\bar\xx^{1-a}\uu^R\|_{L^\infty}^2\left( \|\bar\xx^a\na\n^R\|_{L^p(B_R)}^2 +\|\bar\xx^a\na\hh^R\|_{L^p(B_R)}^2\right)dt\no
\leq &C.\nonumber
\end{align}
Next, one derives from \eqref{3.3} and \eqref{3.24} that
\begin{align}\la{5.16}
&\sup_{0\leq t\leq T_0}t\|\sqrt{\ti\n^R}\ti\uu_t^R\|_{L^2(\mr^2)}^2 +\int_0^{T_0}t\|\na\ti\uu_t^R\|_{L^2(\mr^2)}^2dt\\
\leq &C+C\int_0^{T_0}t\|\na\uu_t^R\|_{L^2(\O)}^2 dt\no
\leq &C.\nonumber
\end{align}

With all these estimates \eqref{5.11}--\eqref{5.16} at hand, we find that the sequence $(\ti\n^R, \uu^R, \ti{\mathbf{\hh}}^R)$ converges, up to the extraction of subsequences, to some limit $(\n, \uu, \hh)$ in the obvious weak sense, that is, as $R\rightarrow\infty$, we have
\be\la{5.17}
\bar\xx\ti{\n^R}\rightarrow \bar\xx\n, \quad \bar\xx\ti{\hh}^R\rightarrow \bar\xx\hh, \quad {\rm in\ }C(\overline{B_N}\times[0, T_0]), \ {\rm for\ any\ } N>0,
\ee
\be\la{5.18}
\bar\xx^a\ti{\n}^R\rightharpoonup \bar\xx^a\n,\quad \bar\xx^a\ti{\hh}^R\rightharpoonup \bar\xx^a\hh, \quad{\rm weakly\ *\ in\ } L^\infty(0, T_0; H^1(\mr^2)\cap W^{1, q}(\mr^2)),
\ee
\be\la{5.19}
\sqrt{\ti\n^R}\ti\uu^R\rightharpoonup\sqrt\n\uu,\quad \na\ti\uu^R\rightharpoonup\na\uu,\quad {\rm weakly\ *\ in\ } L^\infty(0, T_0; L^2(\mr^2)),
\ee
\be\la{5.20}
\na^2\ti\uu^R\rightharpoonup\na^2\uu,\quad{\rm weakly\ in\ } L^{(q+1)/q}(0, T_0; L^q(\mr^2))\cap L^2((0, T_0)\times\mr^2),
\ee
\be\la{5.21}
t^{1/2}\na^2\ti\uu^R\rightharpoonup t^{1/2}\na^2\uu,\quad{\rm weakly\ in\ } L^{2}(0, T_0; L^q(\mr^2)),
\ee
\be\la{5.22}
t^{1/2}\sqrt{\ti\n^R}\ti\uu^R_t\rightharpoonup t^{1/2}\sqrt\n\uu_t, \quad \na\ti\uu^R\rightharpoonup\na\uu,\quad{\rm weakly\ *\ in\ } L^{\infty}(0, T_0; L^2(\mr^2)),
\ee
\be\la{5.23}
t^{1/2}\na\ti\uu^R_t\rightharpoonup t^{1/2}\na\uu_t,\quad{\rm weakly\ in\ } L^{\infty}((0, T_0)\times\mr^2),
\ee
and
\be\la{5.24}
\bar\xx^a\n\in L^\infty(0, T_0; L^1(\mr^2)),\quad\inf_{0\leq t\leq T_0} \int_{B_{2N_0}}\n(t, \xx)d\xx\geq \frac14.
\ee

Then, letting $R\rightarrow\infty$, it follows from \eqref{5.17}--\eqref{5.24} that $(\n, \uu, \hh)$ is a strong solution of \eqref{1.1}--\eqref{1.3} on $(0, T_0]\times\mr^2$ satisfying \eqref{1.6} and \eqref{1.7}. Therefore, we have got the existence part of Theorem \ref{thm1}. However, we  can mimic the argument for the uniqueness of Theorem 1.1 in \cite{lv} to obtain the unique result  and then we finish the proof of Theorem \ref{thm1}.\hfill $\Box$

{\it Proof of Theorem \ref{thm2}.} Let $(\n_0, \uu_0, \hh_0)$ be as in Theorem \ref{thm2}. Without loss of generality, assume that
$$
\int_{\mr^2}\n_0d\xx=1,
$$
which implies that there exists a positive constant $N_0$ such that \eqref{5.1} holds. We construct that $\n_0^R=\hat\n_0^R+R^{-1}e^{-|\xx|^2}$ where $0\leq \hat\n_0^R\in C_0^\infty(\mr^2)$ satisfies \eqref{5.2} and
\be\la{5.33}
\left\{
\begin{array}{lll}
\na^2\hat\n_0^R\rightarrow\na^2\n_0, &{\rm in\ } L^q(\mr^2),\\
\bar\xx^{\delta_0}\na^2\hat\n_0^R\rightarrow \bar\xx^{\delta_0}\na^2\n_0,  &{\rm in\ } L^2(\mr^2),
\end{array}
\right.
\ee
as $R\rightarrow\infty$. Then, we also choose $\hh_0^R\in \{\mathbf{w}\in C_0^\infty(B_R)|\div \mathbf{w}=0\}$ satisfying \eqref{5.3} and
\be\la{5.34}
\left\{
\begin{array}{lll}
\na^2 \hh^R_0\rightarrow\na^2\hh_0,&{\rm in\ } L^q(\mr^2),\\
\bar\xx^{\delta_0}\na^2\hh_0^R\rightarrow \bar\xx^{\delta_0}\na^2\hh_0, &{\rm in\ } L^2(\mr^2),
\end{array}
\right.\quad {\rm as\ } R\rightarrow\infty.
\ee

Then, we consider the unique smooth solution $\uu_0^R$ of the following Stokes problem
\be\la{5.35}
\left\{
\begin{array}{lll}
-\mu\triangle\uu_0^R+\na\pi_0^R=(\nabla\times\hh_0^R)\times\hh_0^R -\n_0^R\uu_0^R+\sqrt{\n_0^R} \mathbf{h}^R, & {\rm in\ } B_R,\\
\div\uu=0,& {\rm in\ } B_R,\\
\uu_0^R=0, &{\rm on\ } \p B_R,
\end{array}
\right.
\ee
where $\mathbf{h}^R=(\sqrt{\n_0}\uu_0+\mathbf{g})*j_{1/R}$ with $j_\delta$ being the standard mollifying kernel of width $\delta$. Multiplying \eqref{5.35} by $\uu_0^R$ and integrating the resultant equation over $B_R$ show that
\begin{align*}
&\|\sqrt{\n_0^R}\uu_0^R\|_{L^2(B_R)}^2+\mu\|\na\uu_0^R\|_{L^2(B_R)}^2 \no
\leq &\int_{\O}|\hh_0^R||\na \hh_0^R||\uu_0^R|d\xx+\|\sqrt{\n_0^R}\uu_0^R\|_{L^2(B_R)} \|\mathbf{h}^R\|_{L^2(B_R)}\no
\leq &\ve\left(\|\sqrt{\n_0^R}\uu_0^R\|_{L^2(B_R)}^2+\|\na\uu_0^R \|_{L^2(B_R)}^2 \right)+C(\ve),
\end{align*}
which implies that
\be\la{5.36}
\|\sqrt{\n_0^R}\uu_0^R\|_{L^2(B_R)}^2+\|\na\uu_0^R \|_{L^2(B_R)}^2 \leq C,
\ee
for some constant $C$ independent of $R$. From \eqref{2.7}, we have
\begin{align}\la{5.37}
&\|\na^2\uu_0^R\|_{L^2(B_R)}+\|\na \pi_0^R\|_{L^2(\O)}\\
\leq &C\||\hh_0^R||\na\hh_0^R|\|_{L^2(\O)}+C\|\n_0^R\uu_0^R\|_{L^2(\O)} +C\|\sqrt{\n_0^R}\mathbf{h}^R \|_{L^2(\O)}\no
\leq &C.\nonumber
\end{align}

Next, extending $\uu_0^R$ to $\mr^2$ by defining $0$ outside $\O$ and denoting it by $\ti\uu_0^R$, we deduce from \eqref{5.36} and \eqref{5.37} that
$$
\|\na\ti\uu_0^R\|_{H^1(\mr^2)}\leq C,
$$
which with \eqref{5.33} and \eqref{5.36} shows that there exists a subsequence $R_j\rightarrow\infty$ and a function $\ti\uu_0\in \{\ti\uu_0\in H^2_{\rm loc}(\mr^2)|\sqrt{\n_0}\ti\uu_0 \in L^2(\mr^2), \na\ti\uu_0\in H^1(\mr^2)\}$ and $\nabla\ti\pi_0\in L^2(\mr^2)$ such that
\be\la{5.38}
\left\{
\begin{array}{lll}
\sqrt{\n_0^{R_j}}\ti\uu_0^{R_j}\rightharpoonup \sqrt{\n_0}\ti\uu_0,&{\rm weakly\ in\ } L^2(\mr^2),\\
\na\ti\uu_0^{R_j}\rightharpoonup \na\ti\uu_0, &{\rm weakly\ in\ } H^1(\mr^2),\\
\nabla\pi^R_0\rightharpoonup \na\ti\pi_0,&{\rm weakly\ in\ } L^2(\mr^2).
\end{array}
\right.
\ee
It is easy to check that $\ti\uu_0^R$ satisfies \eqref{5.35}, then one can deduce from \eqref{5.33}, \eqref{5.34}, \eqref{5.35} and \eqref{5.38} that the pair $(\ti\uu_0,\ti\pi_0)$ satisfies
$$
-\mu\triangle\ti\uu_0+\na\ti\pi_0+\n_0\ti\uu_0=(\nabla\times\hh_0^R)\times\hh_0^R+\n_0\uu_0 +\sqrt{\n_0}\mathbf{g},
$$
which combined with \eqref{1.9} yields that
\be\la{5.39}
\ti\uu_0=\uu_0.
\ee

Next, we get from \eqref{5.35} that
$$
\underset{R_j\rightarrow\infty}{\lim\sup} \int_{\mr^2}\left(|\na \ti\uu_0^{R_j}|^2+\n_0^{R_j}|\ti\uu_0^{R_j}|^2\right)d\xx \leq \int_{\mr^2}\left(|\na\uu_0|^2+\n_0|\uu_0|^2\right)d\xx,
$$
which combined with \eqref{5.38} shows
$$
\lim_{R_j\rightarrow\infty}\int_{\mr^2}|\na \ti\uu_0^{R_j}|^2d\xx= \int_{\mr^2}|\na\uu_0|^2d\xx,\quad\lim_{R_j\rightarrow\infty}\int_{\mr^2} \n_0^{R_j}|\ti\uu_0^{R_j}|^2 d\xx=\int_{\mr^2}\n_0|\uu_0|^2d\xx.
$$
This, along with \eqref{5.38} and \eqref{5.39}, implies  that
\be\la{5.40}
\lim_{R\rightarrow\infty}\left(\|\na(\ti\uu_0^R-\uu_0)\|_{L^2(\mr^2)}+ \|\sqrt{\n_0^R}\ti\uu_0^R-\sqrt{\n_0}\uu_0\|_{L^2(\mr^2)}\right)=0.
\ee
Similar to \eqref{5.40}, we can also obtain that
$$
\lim_{R\rightarrow\infty}\|\na^2(\ti\uu_0^R-\uu_0)\|_{L^2(\mr^2)}=0.
$$

Finally, in terms of Lemma \ref{l2.1}, the initial-boundary value problem \eqref{2.2} with the initial data $(\n_0^R, \uu_0^R, \hh_0^R)$ has a classical solution $(\n^R, \uu^R, \hh^R)$ on $[0, T_R]\times\O$. Hence, there has a generic positive constant $C$ independent of $R$ such that all those estimates stated in Proposition \ref{prop} and Lemma \ref{l4.1}--\ref{l4.4} hold for $(\n^R, \uu^R, \hh^R)$. Extending $(\n^R, \uu^R, \hh^R)$ by zero on $\mr^2\setminus \O$ and denoting
$$
\ti\n^R\triangleq\varphi_R\n^R, \quad \ti\uu^R, \quad \ti\hh^R\triangleq\varphi_R\hh^R,
$$
with $\varphi_R$ as in \eqref{3.6}. We deduce from \eqref{3.3} and Lemma \ref{l4.1}--\ref{l4.4} that the sequence $(\ti\n^R, \ti\uu^R, \ti\hh^R)$ converges weakly, up to the extraction of subsequences, to some limit $(\n, \uu, \hh)$ satisfying \eqref{1.6}, \eqref{1.7} and \eqref{1.10}. Moreover, standard arguments shows that $(\n, \uu, \hh)$ is in fact a classical solution to the problem \eqref{1.1}--\eqref{1.3}. The proof of
Theorem \ref{thm2} is finished. \hfill $\Box$

\section*{Acknowledge}
This work  of M.  Chen is partially supported by the Independent Innovation Foundation of Shandong University (No. 2013ZRQP001) and the National Science Foundation of Shandong province of China under grant (No. ZR2015AM019) and the National Science Foundation of China under grant (No. 11471191). The work of A. Zang is supported in part  National Natural Science Foundation of China (11571279) and a part of Project GJJ151036 supported by Education Department of Jiangxi Province.

\end{document}